\documentclass[12pt]{amsart} 
\usepackage{tikz-cd}
\usepackage{tikz}
\usepackage[utf8]{inputenc}
\usepackage{calculator}

\usepackage{graphicx}

\usepackage{graphicx}
\usepackage{amsfonts}

\usepackage{amsmath}
\usepackage{amssymb}
\usepackage{amsthm}
\usepackage{hyperref}

\tikzcdset{scale cd/.style={every label/.append style={scale=#1},
    cells={nodes={scale=#1}}}}
\numberwithin{equation}{section}
\newcommand{\nn}{\nonumber}
\newcommand{\Perv}{{\rm Perv}}
\newcommand{\MMHM}{{\rm MMHM}}
\newcommand{\MHM}{{\rm MHM}}

\newcommand{\Ind}{{\rm Ind}}

\newcommand{\Id}{{\rm Id}}
\newcommand{\Hom}{{\rm Hom}}
\newcommand{\Crit}{{\rm Crit}}

\newcommand{\ie}{i.e. }
\newlength{\abstractwidth}
\ifx\notheorems\undefined
\theoremstyle{plain}
\newtheorem{theorem}{Theorem}[section]
\newtheorem{proposition}[theorem]{Proposition}
\newtheorem{lemma}[theorem]{Lemma}
\newtheorem{corollary}[theorem]{Corollary}

\theoremstyle{definition}

\newtheorem{definition}[theorem]{Definition}

\theoremstyle{remark}

\title{Hyperbolic localization of the Donaldson-Thomas sheaf}

\author{Pierre Descombes}
\title{Hyperbolic localization of the Donaldson-Thomas sheaf}
\date{}

\begin{document}

\maketitle

\begin{abstract}
    In this paper we prove a toric localization formula in cohomological Donaldson Thomas theory. Consider a $-1$-shifted symplectic algebraic space with a $\mathbb{C}^\ast$ action leaving the $-1$-shifted symplectic form invariant. This includes the moduli space of stable sheaves or complexes of sheaves on a Calabi-Yau threefold with a $\mathbb{C}^\ast$-invariant Calabi-Yau form, or the intersection of two $\mathbb{C}^\ast$-invariant Lagrangians in a symplectic space with a $\mathbb{C}^\ast$-invariant symplectic form. In this case we express the restriction of the Donaldson-Thomas perverse sheaf (or monodromic mixed Hodge module) defined by Joyce et al. to the attracting variety as a sum of cohomological shifts of the DT perverse sheaves on the $\mathbb{C}^\ast$ fixed components. This result can be seen as a $-1$-shifted version of the Białynicki-Birula decomposition for smooth schemes.
\end{abstract}

\tableofcontents

\section{Introduction}

\subsection*{Overview}

A $-1$ shifted symplectic space $X$, such as the moduli space of stable sheaves or complexes of sheaves on a Calabi-Yau threefold or an intersection of two Lagrangians, is informally an algebraic space where the obstruction space is dual to the tangent space. Building on ideas of Kontsevich-Soibelman, Joyce and collaborators have proved that such a space admits a canonical d-critical structure $s$, providing étale charts of $X$ as the critical locus of a function in a smooth space, and defined a perverse sheaf $P_{X,s}$ (and furthermore a monodromic mixed Hodge module on this perverse sheaf) on such $X$ once one has fixed an orientation data. The Euler number of the cohomology of this so-called Donaldson-Thomas sheaf gives the numerical Donaldson-Thomas invariants defined by Thomas. The numerical DT invariant localizes under a torus action, leaving the shifted symplectic structure invariant as the Euler number of a smooth space, i.e. the numerical DT invariant of a space $X$ with a $\mathbb{C}^\ast$ action is the signed sum of the invariants of the $\mathbb{C}^\ast$ fixed components. On a smooth space $U$, the Białynicki-Birula decomposition gives a refinement of the localization formula of the Euler number, which can be expressed using the hyperbolic localization functor as
\begin{align}
        p_!\eta^\ast\mathbb{Q}_U\simeq\bigoplus_{\pi\in\Pi}\mathbb{L}^{d^+_\pi}\mathbb{Q}_{U^0_\pi}
    \end{align}
where $\mathbb{Q}_Y$ is the constant sheaf on $Y$ (or the constant mixed Hodge module), $\mathbb{L}^{1/2}$ is the shift $[-1]$ of a complex (or the underlying monodromic mixed Hodge structure), $\eta$ is the inclusion of the attracting variety $U^+$ of $U$, $p$ is the projection of the attracting variety on the $\mathbb{C}^\ast$ fixed locus with connected component decomposition $U^0:=\bigsqcup_{\pi\in\Pi}U^0_\pi$ and $d^+_\pi$ is the number of contracting weights in the action of $\mathbb{C}^\ast$ on the tangent space of $X$ in $X^0_\pi$. For an algebraic space $X$ with $\mathbb{C}^\star$-invariant d-critical structure $s$ and orientation, the fixed components $X^0_\pi$ have natural d-critical structures $s^0_\pi$ and orientation. In this paper we prove a similar hyperbolic localization formula for the Donaldson-Thomas perverse sheaf (and the monodromic mixed Hodge module) defined using these oriented d-critical structures:
\begin{align}
        p_!\eta^\ast P_{X,s}\simeq\bigoplus_{\pi\in\Pi}\mathbb{L}^{\Ind_\pi/2}P_{X^0_\pi,s^0_\pi}
    \end{align}
where $\Ind_\pi$ is the signed number of contracting weights in the tangent obstruction complex of $X$ on the component $X^0_\pi$. We define the cohomological DT invariant $[X]^{vir}:=[H_c(X,P_{X,s})]$ as the class of hypercohomology with compact support of $P_{X,s}$ in the Grothendieck group of monodromic mixed Hodge structures (in particular, one can extract a virtual Hodge polynomial from it). To compute the cohomological DT invariants of $X$ one is reduced to computing the cohomological DT invariants of the fixed components (which in nice cases are isolated fixed points), and the cohomological DT invariants of the open complement $X-\eta(X^+)$ of the attracting locus (which in nice cases is simpler than the whole moduli space):
\begin{align}\label{nonproj}
    [X]^{vir}=\sum_{\pi\in\Pi}\mathbb{L}^{\Ind_\pi/2}[X^0_\pi]^{vir}+[X-\eta(X^+)]^{vir}
\end{align}
An example where this procedure can be carried out completely is given by framed invariants of toric quivers, as studied in \cite{Descombes2021MotivicDI}. If $X$ is projective, or more generally if each point of $X$ is in the attracting locus of a fixed component (one says that the $\mathbb{C}^\ast$ action is circle-compact), we obtain a formula similar to the localization formula for K-theoretic invariants:
\begin{align}\label{circomploc}
    [X]^{vir}=\sum_{\pi\in\Pi}\mathbb{L}^{\Ind_\pi/2}[X^0_\pi]^{vir}
\end{align}

\subsection*{The Donaldson-Thomas perverse sheaf}

Donaldson-Thomas theory was first developed to count sheaves on Calabi-Yau threefolds. In \cite{thomas1998holomorphic}, Thomas defined the numerical Donaldson-Thomas invariants, giving the virtual Euler number of the moduli space of stable coherent sheaves on a Calabi-Yau threefold, using the perfect obstruction theory given by the Serre duality on the Ext spaces of the sheaves. In \cite{behrend_donaldson-thomas}, Behrend gave a new interpretation of these invariants: expressing the moduli space locally as the critical locus of a potential, the numerical Donaldson-Thomas invariant is given by an Euler characteristic weighted at each point by the Milnor number. In \cite{ks} and \cite{KonSol10}, Kontsevich and Soibelman sketched the definition of a cohomological refinement of this counting, with value in the abelian category of monodromic mixed Hodge modules (MMHM), using the functor of vanishing cycles of this potential, in a partially conjectural framework.

In papers \cite{Joyce2013ACM}, \cite{darbscheme},\cite{Joycesymstab} and \cite{darbstack}, Joyce and collaborators have developed a rigorous cohomological Donaldson-Thomas theory, using the language of $-1$-shifted symplectic structures in derived geometry introduced in \cite{shifsymp}. Informally, derived geometry replaces the notion of tangent space $T_X$ by a tangent complex $\mathbb{T}_X$, where $(\mathbb{T}_X)_0$ gives the tangent directions, $(\mathbb{T}_X)_1$ gives the obstructions, $(\mathbb{T}_X)_2$, $(\mathbb{T}_X)_3$, . ... give the higher obstructions, and $(\mathbb{T}_X)_{-1}$ gives the infinitesimal automorphisms. A $-1$-shifted symplectic structure gives a pairing between $\mathbb{T}_X$ and $\mathbb{T}_X[-1]$: in particular, it pairs tangent directions with obstructions. In \cite{shifsymp} it was shown that moduli stacks of complexes on a Calabi-Yau 3-fold, and intersections of two Lagrangians in a symplectic space, naturally have a $-1$-shifted symplectic structure.

In \cite{darbscheme} it was shown that $-1$-shifted derived algebraic space can be written étale locally as the critical locus of a functional on a smooth scheme, \ie can be covered by critical charts of the form $(R,U,f,i)$, where $R\subset X$ is open, $U$ is a smooth scheme, $f: U\to \mathbb{C}$ is a regular map, and $i:R\to U$ is a closed embedding such that $i(R)=\Crit(f)$, and there is a way to compare intersecting critical charts. More formally, such an algebraic space has a d-critical structure $s$, a notion we will introduce below. In \cite{Joycesymstab}, Joyce and collaborators constructed the Donaldson-Thomas perverse sheaf $P_{X,s}$ carrying a monodromic mixed Hodge module for a d-critical algebraic space $(X,s)$, with additional data called the orientation $K_X^{1/2}$. Its restriction to a critical $(R,U,f,i)$ is given by $\mathcal{PV}_{U,f}$, the perverse sheaf (or monodromic mixed Hodge module) of vanishing cycles defined by applying the vanishing cycles functor $\phi_f$ of $f$ to $\mathcal{IC}_U$, the intersection cohomology complex of $U$ with its natural mixed Hodge module structure, and then restricting to $R$. The orientation is used to glue the perverse sheaves on intersections of critical charts. The numerical Donaldson-Thomas invariant defined in \cite{thomas1998holomorphic} is then, as expected, the Euler number of the cohomology of this perverse sheaf. We will use the formalism of d-critical structures in this article, but not the formalism of derived geometry and shifted symplectic structures.

\subsection*{Hyperbolic localization}
The aim of this paper is to provide a way to compute the cohomological Donaldson-Thomas invariants by localization, namely, given a d-critical algebraic space $X$ with a $\mathbb{C}^\ast$ action and a $\mathbb{C}^\ast$ invariant d-critical structure, express the invariants of $X$ in terms of the invariants of the $\mathbb{C}^\ast$ fixed algebraic space $X^0$. Graber and Pandharipande proved a torus localization formula for numerical Donaldson-Thomas invariants in \cite{Graber1997LocalizationOV}. An analogue formula was derived in \cite{Behrend08symmetricobstruction} using the alternative definition with weighted Euler characteristic: the numerical Donaldson-Thomas invariant of $X$ is the sum of those of each component of $X^0$ weighted by a sign given by the parity of the dimension of the normal space to the component. Thus, numerical Donaldson-Thomas invariants localize under torus action exactly like the Euler numbers of smooth spaces.

In \cite{BiaynickiBirula1973SomeTO}, Białynicki-Birula proved that, on a smooth scheme $X$ with $\mathbb{C}^\ast$ action, the attracting subset of $X$ (\i.e. the subset of points $x\in X$ such that $\lim_{t\to 0}t.x$ exists) admits a kind of cellular decomposition. Each cell of this so-called Białynicki-Birula decomposition is an affine fiber bundle on a component of the fixed variety $X^0$, whose dimension is given by the number of contracting $\mathbb{C}^\ast$-weights in the tangent space of $X^0$. Thus one can compute the cohomology of the attracting subset of $X$ in terms of the cohomology of $X^0$. In \cite{Braden2002HyperbolicLO} and \cite{drinfeld2013algebraic}, Braden and Drinfeld gave a functorial reformulation of the Białynicki-Birula decomposition that can be extended to singular algebraic spaces. For any algebraic space $X$ with $\mathbb{C}^\ast$ action one can define the $\mathbb{C}^\ast$ fixed space $X^0$ and also the contracting (resp. repelling) space $X^\pm$ in an obvious way using the functor of points. Then there is the so-called hyperbolic localization diagram:
\[\begin{tikzcd}
X & X^\pm\arrow[l,swap,"\eta^\pm_X"]\arrow[r,"p^\pm_X"] & X^0
\end{tikzcd}\]
where $\eta^\pm_X$ is the inclusion of the attracting (resp. repelling) subset, and $p^+_X$ (or $p^-_X$) sends $x$ to $\lim_{t\to 0}x$ (or $\lim_{t\to \infty}x$). This diagram is easy to describe if $X$ is an affine scheme: it suffices to find a $\mathbb{C}^\ast$-equivariant covering of $X$ by affine spaces. The existence of such a covering is not at all obvious, and is far too restrictive in Zariski topology. Fortunately, a deep result of \cite{Alper2015AL} and \cite{Alper2019TheL} shows that quasi-separated algebraic spaces of finite type locally have a $\mathbb{C}^\ast$-equivariant covering by affine spaces. Under these assumptions $X^0$ and $X^\pm$ are then algebraic spaces. The functors $(p^\pm_X)_!(\eta^\pm_X)^\ast$ are called hyperbolic localization functors. From \cite{drinfeld2013algebraic} there is a natural morphism built from morphisms of the six-functor formalism in the derived category of constructible sheaves, or the derived category of monodromic mixed Hodge modules:
\begin{align}
    S_X:\mathbb{D}(p^-_X)_!(\eta^-_X)^\ast\mathbb{D}\to (p^+_X)_!(\eta^+_X)^\ast
\end{align}
This is by \cite{drinfeld2013algebraic} and \cite{Ric16} an isomorphism when these functors are applied to $\mathbb{C}^\ast$-equivariant complexes, \ie complexes equivariant under the action of $\mathbb{C}^\ast$, in a sense we will explain below.\medskip

Assume now that $X$ is smooth, and denote by $\Pi$ the cocharacters of $\mathbb{C}^\ast$, and for $\pi\in\Pi$ $X^0_\pi$ the union of the connected components of $X^0$, where the $\mathbb{C}^\ast$ action on $T_{X^0}$ has the cocharacter $\pi$. For $\pi\in\Pi$, denote by $d^0_\pi,d^+_\pi,d^-_\pi$ the number of invariants, resp. contracting, resp. repelling weights in $\pi$, and define $\Ind_\pi:=p^+_\pi-p^-_\pi$. The Białynicki-Birula decomposition for smooth schemes says that $X^\pm$ is a disjoint union of affine fiber bundles $X^\pm_\pi$ of dimension $d^\pm_\pi$ over $X^0_\pi$. Using the hyperbolic localization functors, this can be reformulated as
\begin{align}
    (p^\pm_X)_!(\eta^\pm_X)^\ast\mathbb{Q}_U&=\bigoplus_{\pi\in\Pi}\mathbb{L}^{d^\pm_\pi}\mathbb{Q}_{U^0_\pi}\nn\\
    \Rightarrow (p^\pm_X)_!(\eta^\pm_X)^\ast\mathcal{IC}_U&=\bigoplus_{\pi\in\Pi}\mathbb{L}^{\pm\Ind_\pi/2}\mathcal{IC}_{U^0_\pi}\label{locintcoh}
\end{align}
where $\mathbb{Q}_Y$ and $\mathcal{IC}_Y$ are respectively the constant sheaf and the intersection complex of $Y$, or their natural mixed Hodge modules. This isomorphism sends $S_X$ to the self-duality isomorphism of the intersection complex.

\subsection*{The main result}
For a $\mathbb{C}^\ast$-equivariant oriented d-critical algebraic space $X$ satisfying the assumptions under which the hyperbolic localization functors were defined, we will prove that $X^0$ has a natural d-critical structure $s^0$ and a natural orientation $K_{X^0}^{1/2}$ on $(X^0,s^0)$ induced by $K_X^{1/2}$.  The aim of this article is to prove that an analogue of \eqref{locintcoh} holds for the Donaldson-Thomas perverse sheaf (res. monodromic mixed Hodge module) on a d-critical oriented algebraic space :\medskip

\begin{theorem}(Theorem \ref{theospace})
    For $X$ a quasi-separated locally of finite type algebraic space with an action of a one-dimensional torus $\mathbb{C}^\ast$ and $\mathbb{C}^\ast$-equivariant d-critical structure $s$ and orientation, there are natural isomorphisms of perverse sheaves and monodromic mixed Hodge modules:
    \begin{align}
        \beta^\pm_{X,s}:(p_X^\pm)_!(\eta_X^\pm)^\ast P_{X,s}\overset{\simeq}{\to}\bigoplus_{\pi\in\Pi}\mathbb{L}^{\pm \Ind_\pi/2}P_{X^0_\pi,s^0_\pi}
    \end{align}
    where $\Ind_\pi$ is the signed number of contracting weights in the $\mathbb{C}^\ast$ action on the tangent obstruction complex of $X$ at $X^0_\pi$.
\end{theorem}

\subsection*{Sketch of the proof}

The main ingredient of the proof of this theorem is the commutation of the hyperbolic localization and vanishing cycles functor proved in \cite{Ric16}. On a critical chart $(R,U,f,i)$, denoting by $(R^0,U^0,f^0,i^0)$ the fixed critical chart, combining with the classical Białynicki-Birula decomposition of $U$ \eqref{locintcoh}, one obtains isomorphisms in the derived category of constructible sheaves with monodromy, or of monodromic mixed Hodge modules:
\begin{align}
        \beta^\pm_{U,f}:(p_R^\pm)_!(\eta_R^\pm)^\ast \mathcal{PV}_{U,f}\to\bigoplus_{\pi\in\Pi}\mathbb{L}^{\pm \Ind_\pi/2}\mathcal{PV}_{U^0_\pi,f^0_\pi}
\end{align}
This is the content of Proposition \ref{propchi}. We need to show that these isomorphisms behave well with respect to pullback by étale maps of the critical locus: this is the content of Proposition \ref{propcomsmooth}. We also need to show that they behave well with respect to the Thom-Sebastiani isomorphism of \cite{MasThomSeb}:
\begin{align}
    \mathcal{TS}_{U,f,V,g}:\mathcal{PV}_{U\times V,f\boxplus g}\simeq\mathcal{PV}_{U,f}\boxtimes\mathcal{PV}_{S,g}
\end{align}
This is the content of Proposition \ref{propcomthomseb}.

Now consider a d-critical oriented algebraic space $(X,s)$. The isomorphism $\beta^\pm_{X,s}$ is an isomorphism in the abelian category $\bigoplus_\pi \Perv(X^0_\pi)[-\Ind_\pi]$ of shifted perverse sheaves, or the abelian category $\bigoplus_\pi \MMHM(X^0_\pi)[-\Ind_\pi]$ of monodromic mixed Hodge modules, which form stacks for the étale topology on $X$, so it suffices to define them as $\beta^\pm_{U,f}$ on any critical chart $(R,U,f,i)$, and to show that they agree on intersections of critical charts, compatibility with monodromy and self-duality can be checked on critical charts. 
Note that the tangent obstruction complex of $X$ on $R$ is quasi-isomorphic with $0\to T_U\to T^\ast_U\to 0$, so $\Ind_\pi$ is the signed number of contracting weights in the tangent-obstruction complex of $X$ at $X^0_\pi$.\medskip

For d-critical algebraic spaces, Joyce proved in \cite{Joyce2013ACM} that two critical charts intersecting at a point can be related in an étale neighbourhood of that point by stabilization by quadratic forms, that is, by embeddings $(U,f)\hookrightarrow (U\times E,f\boxplus q)$, where $q$ is a non-degenerate $\mathbb{C}^\ast$-invariant quadratic form on a $\mathbb{C}^\ast$-equivariant vector space $E$. In \cite{Joycesymstab}, Joyce et al. glue the perverse sheaves of vanishing cycles using the isomorphism:
\[\begin{tikzcd}
        \mathcal{PV}_{U\times E,f\boxplus q}\arrow[r,"\mathcal{TS}_{U,f,E,q}"] & \mathcal{PV}_{U,f}\boxtimes\mathcal{PV}_{E,q}\arrow[r,"\simeq"] & \mathcal{PV}_{U,f}
\end{tikzcd}\]
where they use the isomorphism $\mathcal{PV}_{E,q}\simeq\mathbb{Q}$ with a choice of orientation on $E$ (this choice is at the root of orientation problems in cohomological DT theory). To prove that the isomorphisms $\beta^\pm$ agree on the intersection of critical charts, we prove that they are compatible with stabilization, which is done using compatibility with the Thom-Sebastiani isomorphism and a careful comparison of orientations.\medskip

In this paper we use the functorial properties of the six-functor formalism in the derived categories $D^b_c(Y)$ of bounded constructible complexes over various algebraic spaces $Y$, or, more formally, a $(\infty,1)$-categorical enhancement of them. When we say that a diagram commutes from the naturality of the six-functor formalism, we imply that it commutes up to a natural 2-isomorphism. Since the final diagrams lie in the abelian heart $\bigoplus_\pi \Perv(\mathcal{X}^0_\pi)[-\Ind_\pi]$ of $D^b_c(X^0)$, which is a classical category, they really commute. Moreover, it is assumed without specification that all algebraic spaces are quasi-separated and locally of finite type.\medskip

For $X$ an algebraic space, the category $\MHM(X)$ of mixed Hodge modules introduced in \cite{Saito1990MixedHM} is an abelian category with a faithful and exact functor $rat:\MHM(X)\to \Perv(X)$. In particular, because $rat$ is faithful, the commutativity of a diagram can be checked at the level of perverse sheaves. Moreover, since a mixed Hodge module $M$ is the data of filtrations on the $\mathcal{D}$ module associated with the complexification of the perverse sheaf $rat(M)$ by the Riemann-Hilbert correspondence, a morphism of mixed Hodge modules is an isomorphism if the underlying morphism of perverse sheaves is an isomorphism. There is a six-functor formalism on the derived category of mixed Hodge modules which, under $rat$, projects to the six-functor formalism on $D^b_c(X)$, the derived category of constructible complexes. Furthermore, since $rat$ is exact, a morphism between complexes of mixed Hodge modules is an isomorphism in the derived category if the underlying morphism of constructible complexes is an isomorphism in the derived category. Thus, the extension of our results to monodromic mixed Hodge modules is straightforward: our morphisms are easily defined using the six-functor formalism, so we used the same definition to extend them to morphisms in the derived category of monodromic mixed Hodge modules. The hard work is to check that these morphisms are isomorphisms or that they give commutative diagrams, but this can be checked at the level of perverse sheaves as explained above.\medskip

\subsection*{Analytic version}

In \cite{Joycesymstab} there is also an analytic version of cohomological Donaldson-Thomas theory, for analytic oriented d-critical analytic spaces which are locally the critical locus of a holomorphic function. Unfortunately, we lack crucial results on $\mathbb{C}^\ast$ actions on analytic spaces. If an analytic space with $\mathbb{C}^\ast$ action has a $\mathbb{C}^\ast$ equivariant analytic (resp. smooth) cover by an affine analytic space, we show that the fixed and attracting/repelling varieties are analytic spaces, hence the hyperbolic localization functors and the morphism $S_X$ are defined. However, we do not know under which generality such a $\mathbb{C}^\ast$-equivariant cover exists. Moreover, we need an analytic version of \cite[Theo B]{Ric16} which says that $S_X$ is an isomorphism when applied to $\mathbb{C}^\ast$-equivariant complexes. Furthermore, according to \cite[Remark 2.20]{Joycesymstab}, the six-functor formalism is not fully defined for mixed Hodge modules in the analytic case, so it might be difficult to lift a localization result from the level of perverse sheaves to the level of mixed Hodge modules. Once such results are available, all the results of this article should extend directly to the analytic case.

\subsection*{Case with semistable objects}\label{sectsemis}

Here we consider that we have $-1$-shifted symplectic spaces, which is the case for the moduli space of stable objects in CY3 categories as proved in \cite{Brav2018RelativeCS}. If one has semistable objects that are not stable, the stack $\mathfrak{X}$ of semistable objects is $-1$-shifted symplectic from \cite{Brav2018RelativeCS}. Joyce and collaborators then built a perverse sheaf $P_{\mathfrak{X}}$ on the $-1$-shifted symplectic stack with orientation data in \cite{darbstack}. If one had a nice six-functor formalism for mixed Hodge modules on stack, one could similarly define a DT mixed Hodge module on any $-1$-shifted symplectic stack, but such a theory is still lacking.

Under very general assumptions, in this setting there is a good moduli space $JH:\mathfrak{X}\to X$ from \cite{Alper2018ExistenceOM} which is not $-1$-shifted symplectic if there are semistable objects which are not stable. As suggested by Toda in \cite{Toda2017GopakumarVafaIA}, one can mimic the definition of \cite{DavMein} to define a BPS perverse sheaf:
\begin{align}
    \mathcal{BPS}:=^p\mathcal{H}^1(JH_!P_{\mathfrak{X}})
\end{align}
where $^p\mathcal{H}^1$ denotes the first perverse cohomology. However, a nice integrality formalism generalising the results of \cite{DavMein} is still missing.

It seems plausible that, once such a formalism is established, a similar localization formula holds for the BPS perverse sheaf of mixed Hodge modules, extending localization techniques to the case with semistable objects, and this is the subject of future work. A localization result for the perverse sheaf $P_{\mathfrak{X}}$ on the stack was proved in the first version of this paper, but we have removed it for the clarity of the present paper. We have only kept the subsection \ref{sectsmooth} proving the compatibility of the localization isomorphism with smooth maps of d-critical loci, which was instrumental in proving the localization formula for stacks, because it may be useful in the future, and do not use any stack formalism.

\subsection*{Relations to other works}\label{relwork}

The idea of using hyperbolic localization to obtain a localization formula in cohomological DT theory was first formulated by Balazs Szendroi in \cite[Section 8.4]{Szendroi:2015xna}. It was applied in \cite[Section 6]{Nakajima2016LecturesOP} and in \cite[Section 8.3]{RSYZ19}, where it was used in a specific example of framed representations of quivers with potential. In \cite{Ric16} Timo Richarz proved the commutativity of hyperbolic localization with the vanishing cycles functor in a more general way.  We used this result to establish the formula \eqref{nonproj} for any critical locus of a potential in \cite{Descombes2021MotivicDI}, and the extension of this result to $-1$-shifted symplectic spaces was suggested to us by Richard Thomas.

In \cite{Bussi2019OnMV} the authors defined a Donaldson-Thomas motive which glues the motive of vanishing cycles defined by arc spaces. Motives glue in Zariski topology and not in étale topology, so to prove a localization formula one must a priori restrict to d-critical schemes admitting Zariski $\mathbb{C}^\ast$-equivariant critical charts. According to \cite[Prop 2. 43]{Joyce2013ACM}, the existence of such charts is equivalent to the existence of a $\mathbb{C}^\ast$-equivariant Zariski cover (such a $\mathbb{C}^\ast$ action is then called 'good'), which holds for normal spaces from Sumihiro's theorem \cite{sumihiro1974equivariant}, but not in greater generality (see \cite[Example 2. 46]{Joyce2013ACM}). Davesh Maulik has proved a formula similar to \eqref{circomploc} for motivic DT invariants on d-critical schemes with a good circle-compact $\mathbb{C}^\ast$ action, as explained in \cite[section 5.3]{darbstack}, in an unpublished preprint (private communication). A generalisation of this result for non-Archimedean geometry was subsequently proved in \cite[Theo 7.17]{jiang2017moduli}. For non-circle compact actions one can compute the DT motive of the attracting variety, which is by definition circle compact, in this way.

A toric localization similar to \eqref{circomploc} also exists for K-theoretic DT invariants as defined in \cite{NekOk}. The K-theoretic DT invariants are a refinement of the numerical DT invariants defined for projective moduli spaces with symmetric obstruction theories, developed in parallel to the motivic and cohomological refinement of Kontsevich-Soibelman and Joyce and collaborators. In general, they are expected to correspond to the $\chi_y$ genus of the Hodge polynomial of cohomological DT theory, so in particular one replaces $\mathbb{L}^{1/2}$ by $-y$ in K-theoretic formulae. If the moduli space $X$ is not projective, but has a $\mathbb{C}^\ast$ action with projective fixed components, it was suggested in \cite{NekOk} to use the equation \eqref{circomploc} to define the $K$-theoretic invariants of $X$. However, this definition depends on the choice of the $\mathbb{C}^\ast$ action (this choice is called the choice of the slope). The equation \eqref{nonproj} in the non-projective case explains the origin of this ambiguity: one computes by toric localization only the virtual cohomology of the attracting variety, which is not the whole moduli space and depends on the chosen $\mathbb{C}^\ast$-action.\medskip

This dependence on the slope was studied explicitly in \cite{Arb} for the moduli space of framed representations of a toric quiver, and this was related to the ambiguity in the refined topological vertex of \cite{Iqbal:2007ii}. In this case there is a two-dimensional torus invariant acting on the moduli space of framed representations, scaling the arrows of the quiver by leaving the potential invariant, so the space of slopes is $\mathbb{P}^1_\mathbb{R}$. The fixed points can be described as molten crystals from \cite{Mozgovoy:2008fd}. In \cite[Prop 3.3]{Arb} it was established that there is a wall and chamber structure on the space of slopes, with the generating functions of framed invariants being constant in a chamber and jumping at a wall, the wall corresponding to slopes where the weight of an elementary cycle of the quiver becomes attracting or repelling. This is rather strange at first sight, because inside a given chamber the cohomological weight of a given molten crystal changes at many walls, but the final result does not change: these walls are 'invisible'. In \cite{Descombes2021MotivicDI} it was established that the attracting variety is the subspace of representations where the cycles with repelling weights are nilpotent, so the attracting variety changes exactly on the walls defined in \cite{Arb}, \ie \eqref{nonproj} give an explanation of this wall and chamber structure. Moreover, using a nilpotent/invertible decomposition for the unframed representation and a wall-crossing relation between framed and unframed invariants, \cite{Descombes2021MotivicDI} obtained the full framed generating series by multiplying the one obtained by localization by a generating series of framed invariants where some cycles are imposed to be invertible. The latter is easy to compute and has a universal closed formula for all toric quivers. Note that in this case the moduli space is the critical locus of the potential of the quiver, so one does not need all the subtleties of gluing, \ie one only needs the proposition \ref{propchi}.

\subsection*{Acknowledgements}
I thank Boris Pioline and Olivier Schiffmann, my PhD supervisors, and David Rydh, Benjamin Hennion, Davesh Maulik, Timo Richarz, Matthieu Romagny, Balazs Szendroi and Richard Thomas for interesting discussions and comments on a preliminary version of the paper.

\section{Classical Białynicki-Birula decomposition}

\subsection{Hyperbolic localization}

Consider an algebraic space $S$, and $X$ be an algebraic space over $S$ with a relative $\mathbb{C}^\ast$ over $S$.\medskip

\begin{definition}
    An algebraic space $X$ over an algebraic space $S$ with a relative $\mathbb{C}^\ast$ over $S$ is said to be étale locally (resp. analytical locally) linearizable if there is a $\mathbb{C}^\ast$-equivariant étale (resp. analytical) covering family $\{U_i\to X\}_i$ where the $U_i$ are affine $S$-algebraic spaces with $\mathbb{C}^\ast$-action.
\end{definition}

According to \cite[Cor 20.2]{Alper2019TheL}, an algebraic space $X$ is then étale locally linearizable over $S$ when $X$ and $S$ are quasi-separated, and $X$ is locally of finite type over $S$, and we expect similarly that the fact of being analytic locally linearizable is not too restrictive. We will here adapt the setting of \cite{Ric16} to the complex analytic case. Consider the following functors on the category of $S$-algebraic spaces:
\begin{align}
    X^0:&Y\mapsto \Hom_S^T(Y,X)\nn\\
    X^+:&Y\mapsto \Hom_S^T((\mathcal{A}^1_Y)^+,X)\nn\\
    X^-:&Y
    \mapsto \Hom_S^T((\mathcal{A}^1_Y)^-,X)
\end{align}
where the superscript $\mathbb{C}^\ast$ denotes the $\mathbb{C}^\ast$-equivariant morphism, and $T,(\mathcal{A}^1_Y)^+,(\mathcal{A}^1_Y)^-$ has the trivial, resp. usual, resp. opposite $\mathbb{C}^\ast$-action. 

The structure morphism $(\mathcal{A}^1_S)^\pm\to S$ is $\mathbb{C}^\ast$-equivariant and defines then a morphism:
\begin{align}
    \zeta^\pm:X^0\to X^\pm
\end{align}
The zero section of $(\mathcal{A}^1_S)^\pm\to S$ defines a morphism:
\begin{align}
    p^\pm:X^\pm\to X^0
\end{align}
such that $p^\pm\circ\zeta^\pm=\Id_{X^0}$. Finally, the unit section of $(\mathcal{A}^1_S)^\pm\to S$ defines a morphism:
\begin{align}
    \eta^\pm:X^\pm\to X
\end{align}
such that $\xi:=\eta^+\circ\zeta^+=\eta^-\circ\zeta^-$ is the inclusion of the subfunctor $X^0\subset X$. Hence it defines the hyperbolic localization diagram as in \cite{Ric16}:
\[\begin{tikzcd}
    X^0 &&&\\
    & X^+\arrow[dl,"\eta^+"]\arrow[ul,"p^+"] &\\
    X && X^+\times_{X_0}X^-\arrow[lu,"'\eta^-"]\arrow[ld,swap,"'\eta^+"] & X^0\arrow[ull,"\zeta^+"]\arrow[dll,swap,"\zeta^-"]\arrow[l,"j"]\\
    & X^-\arrow[ul,swap,"\eta^-"]\arrow[dl,"p^-"] &\\
    X^0&&&
\end{tikzcd}\]

\begin{proposition}
    For $X$ an $S$ algebraic space with an étale (resp. analytical)-locally linearizable $\mathbb{C}^\ast$-action, $X^0$ is representable by a closed algebraic (resp. analytic) subspace of $X$, and $X^\pm$ are representable as $X^0$-affine algebraic (resp. algebraic) spaces, and $j$ is open and closed.
\end{proposition}

Proof: The version for algebraic spaces is \cite[Prop 1.17]{Ric16}. We give a similar proof in the complex analytic space. The proof of the proposition when $X$ is affine goes exactly the same as in \cite[Lem 1.9]{Ric16}, namely $X^0,X^+,X^-$ are represented by the closed subscheme defined by the ideal of homogeneous elements of degree zero (resp. strictly negative, resp. strictly positive positive degree), and $X^0=X^+\times_{X_0}X^-$.\medskip

Suppose now the $U\to X$ is the embedding of a $\mathbb{C}^\ast$-invariant open subset of $X$ in the analytic topology. One can adapt \cite[Lem 1.10]{Ric16}: for $Y/S$ a complex analytic space, an element $\phi\in(U\times_X X_0)(Y)$ corresponds to a commutative diagram:
\[\begin{tikzcd}
    & U\arrow[d]\\
    Y\arrow[ur,"\tilde{f}"]\arrow[r,"f"] & X
\end{tikzcd}\]
and, because $U\to X$ is an open immersion, $\tilde{f}$ is $\mathbb{C}^\ast$-equivariant, hence $\phi\in U^0(T)$. Then $U^0=U\times_X X^0$ a functors, and then, taking a $\mathbb{C}^\ast$-equivariant analytic covering family $\{U_i\to X\}_i$ where the $U_i$ are $S$-affine, one obtains a Cartesian commutative diagram:
\[\begin{tikzcd}
    \bigsqcup_i U^0_i\arrow[r]\arrow[d] & X^0\arrow[d]\\
    \bigsqcup_i U_i\arrow[r] & X
\end{tikzcd}\]
Because the top left element is representable, and the vertical left arrow is a closed immersion, one obtains that $X^0$ is representable as a close analytic subspace of $X$.\medskip

For $U\to X$ the embedding of a $\mathbb{C}^\ast$-invariant open subset of $X$ in the analytic topology, one can adapt \cite[Lem 1.11]{Ric16}: for $Y/S$ a complex analytic space, an element $\phi\in(U^0\times_{X_0} X^\pm)(Y)$ corresponds to a commutative diagram of $\mathbb{C}^\ast$-equivariant morphisms:
\[\begin{tikzcd}
    Y\arrow[r,"f|_Y"]\arrow[d] & U\arrow[d]\\
    (\mathbb{A}^1_Y)^\pm\arrow[r,"f"] & X
\end{tikzcd}\]
where $Y\to(\mathbb{A}^1_Y)^\pm$ is the zero section. The open subset $f^{-1}(U)$ contains $Y$ and is $\mathbb{C}^\ast$ invariant in $\mathbb{A}^1_Y$, hence is the whole $\mathbb{A}^1_Y$, \ie $\phi\in U^\pm(Y)$. Hence $U^\pm=U^0\times_{X_0} X^\pm$ as functors, and taking a $\mathbb{C}^\ast$-equivariant analytic covering family $\{U_i\to X\}_i$ where the $U_i$ are $S$-affine, one obtains a Cartesian commutative diagram:
\[\begin{tikzcd}
    \bigsqcup_i U^\pm_i\arrow[r]\arrow[d] & X^\pm\arrow[d]\\
    \bigsqcup_i U_i^0\arrow[r] & X^0
\end{tikzcd}\]
Because the top left element is representable, and the vertical left arrow is affine, one obtains that $X^0$ is representable as an affine analytic space on $X^0$. Moreover, $X^0=X^+\times_{X_0}X^-$ is an open and closed immersion because it is the case analytical locally.$\Box$\medskip

There is a natural transformation $\mathbb {D}(\eta^-)_!(p^-)^\ast\mathbb{D}\simeq(\eta^-)_\ast(p^-)^!\to(\eta^+)_!(p^+)^\ast$ built using the six-functor formalism between derived functors $D^b_c(X)\to D^b_c(X^0)$, as explained in \cite[Section 2.2]{Ric16}, and we can define the same similar transformation at the level of mixed Hodge modules. We change slightly the terminology of \cite{Ric16}, replacing the term 'monodromic' by 'equivariant' to avoid future confusions with monodromic mixed Hodge modules. Let $a,p:T\otimes_S X \to X$ denote the action (resp. projection). We say that a constructible complex $A\in D^b_c(X)$ is naively $\mathbb{C}^\ast$-equivariant if there exists an isomorphism $a^\ast A\simeq p^\ast A$ in $D^b_c(T\otimes_S)$. Let us define $D^b_c(X)^{\mathbb{C}^\ast-eq}$ to be the full subcategory strongly generated by naively $\mathbb{C}^\ast$-equivariant complexes, i.e. generated by a finite iteration of taking the cone of a morphism in $D^b_c(X)$. The objects in $D^b_c(X)^{\mathbb{C}^\ast-eq}$ are called $\mathbb{C}^\ast$-equivariant. The result \cite[Theo B]{Ric16} can then be rephrased as:\medskip

\begin{proposition}\label{hyplocdual}
    For $X$ an algebraic space with $\mathbb{C}^\ast$ action, the natural morphism of constructible complexes (or complexes of mixed Hodge modules) $S_X(A):\mathbb {D}(\eta^-_X)_!(p^-_X)^\ast\mathbb{D}(A)\to(\eta^+_X)_!(p^+_X)^\ast(A)$ is an isomorphism when $A$ is a $\mathbb{C}^\ast$-equivariant constructible complex.
\end{proposition}

where we have used the fact that a morphism of complexes of mixed Hodge modules is an isomorphism if the underlying morphism of constructible is an isomorphism.

\subsection{Compatibility with smooth morphisms}

The hyperbolic localization diagram is functorial, namely for $\phi:X\to Y$ there are $\mathbb{C}^\ast$-equivariant morphisms $\phi^\pm:X^\pm\to Y^\pm$ and morphisms $\phi^0:X^0\to Y^0$, obtained by composition with $\phi$ using the functor description, such that all the squares commutes in the following diagram:
\[\begin{tikzcd}
    X\arrow[d,"\phi"] & X^\pm\arrow[l,"\eta_X^\pm"]\arrow[r,shift left=1ex,"p_X^\pm"]\arrow[d,"\phi^\pm"] & X^0\arrow[l,"\zeta_X^\pm"]\arrow[d,"\phi^0"]\\
    Y & Y^\pm\arrow[l,"\eta_Y^\pm"]\arrow[r,shift left=1ex,"p_Y^\pm"] & Y^0\arrow[l,"\zeta_Y^\pm"]
\end{tikzcd}\]

\begin{proposition}\label{prophyplocsmoothaff}
    For a smooth $\mathbb{C}^\ast$-equivariant morphism $\phi:X\to Y$, $\phi^\pm$ and $\phi^0$ are also smooth, and the natural map $X^\pm\to X^0\times_{Y^0} Y^\pm$ is an affine fibre bundle.
\end{proposition}

Proof:  Consider first that $\phi$ is étale (resp. the embedding of an open subset in the analytic topology): then \cite[Lem 1.10, 1.11]{Ric16} (resp. our proof of the representability of $X^0,X^\pm$) shows respectively that $X^0=X\times_Y Y^0$ and that $X^\pm=X^0\times_{Y^0} Y^\pm$: in particular, $\phi^0$, and therefore $\phi^\pm$, are étale (resp. the embedding of an open subset in the analytic topology) by base change.\medskip

Consider now the case of an affine $\mathbb{C}^\ast$-equivariant fibration $\phi:Y\times_S \mathbb{A}^d_S\to Y$ where the fiber has a linear $\mathbb{C}^\ast$ action. Then $\phi^0:Y^0\times_S (\mathbb{A}^d_S)^0\to Y^0$ and $\phi^\pm:Y^\pm\times_S (\mathbb{A}^d_S)^\pm\to Y^\pm$ are smooth, and:
\begin{align}
    &(Y\times_S \mathbb{A}^d_S)^\pm=Y^\pm\times_S (\mathbb{A}^d_S)^\pm\nn\\
    \to&  (Y\times_S \mathbb{A}^d_S)^0\times_{Y^0} Y^\pm=Y^\pm\times_S (\mathbb{A}^d_S)^0
\end{align}
is affine.\medskip

Suppose now that $\phi$ is smooth. For $x\in X^0$, we can restricts to a $\mathbb{C}^\ast$-invariant étale (resp. open in the analytic topology) neighborhood of $\phi(x)$ in which $Y$ is affine, and further restricts to a $\mathbb{C}^\ast$-invariant étale (resp. open in the analytic topology) neighborhood of $x$ in $X$ which is affine. As argued in the proof of proposition 2.43 of \cite[Prop 2.43]{Joyce2013ACM}, because $\mathbb{C}^\ast$ is a torus, an affine $\mathbb{C}^\ast$-equivariant space can be written as $Spec(R[x_1,...,x_n]/(h_1,...,h_r))$ where the $x_i$ are $\mathbb{C}^\ast$-equivariant coordinates, and the $h_i$ are $\mathbb{C}^\ast$-equivariant polynomial. Then the $\mathbb{C}^\ast$-equivariant smooth map $\phi:X\to Y$ is étale locally given by the dual of the map of rings:
\begin{align}
    \tilde{R}[x_1,...,x_n]/(h_1,...,h_r)\to \tilde{R}[y_1,...,y_m]/(k_1,...,k_s)
\end{align}
where we replace the ring of polynomial in $n$ variables by the ring of holomorphic functions in $n$ variables i the analytic setting, which can be rewritten:
\begin{align}
    R \to R[y_1,...,y_m]/(l_1,...,l_{r+s})
\end{align}
where $R:=\tilde{R}[x_1,...,x_n]/(h_1,...,h_r)$ and the $l_i$ are the $\mathbb{C}^\ast$-equivariant polynomials (resp. holomorphic functions) $k_1,...,k_r,x_1-\phi'(x_1),...,x_n-\phi'(x_n)$. Denote by $c$ the rank of the matrix $(\frac{\partial l_i}{\partial y_j}|_x)_{ij}$ in the neighborhood of $x$, and up to reordering suppose that $(\frac{\partial l_i}{\partial y_j}|_x)_{0\leq i,j\leq c}$ is invertible at $x$, hence in the neighborhood of $x$ it is still invertible, \ie $R[y_{c+1},...,y_m]\to R[y_1,...,y_m]/(l_1,...,l_{r+s})$ is étale (resp. a local homeomorphism) and $R[y_1,...,y_m]/(l_1,...,l_c)\to R[y_1,...,y_m]/(l_1,...,l_{r+s})$ is a local homeomorphism, hence $\phi$ can be étale locally written as a composition of $\mathbb{C}^\ast$-equivariant morphisms:
\begin{align}
    R\to R[y_{c+1},...,y_m]\to S[y_1,...,y_m]/(l_1,...,l_{r+s})
\end{align}
where the second map is étale (resp. a local homeomorphism). Because $X^0\to Y^0$, $X^\pm\to Y^\pm$ and $X^\pm\to X^0\times_{Y^0} Y^\pm$ are smooth for $X\to Y$ $\mathbb{C}^\ast$-equivariant and étale (resp. open embeddings in the analytic topology) and for $X\to Y$ affine fibrations with $\mathbb{C}^\ast$-linear fibers, then it is true for smooth maps because the property of being smooth can be checked étale locally.\medskip

In order to show that the smooth morphism $q:X^\pm\to X^0\times_{Y^0} Y^\pm$ is an affine fibre bundle, one can take $X=X^\pm$, $Y=Y^\pm$. We provide here a relative version of the proof of \cite[Lem 7.2]{JELISIEJEW2019290}. Consider the sheaf of ring $\mathcal{B}_0$ of $Y^0$ (resp. $\mathcal{C}_0$ of $X^0$). Because $Y^\pm\to Y$ (resp. $X^\pm\to X$) is affine, it corresponds to a sheaf of $\mathbb{Z}$-graded $\mathcal{B}_0$- (resp. $\mathcal{C}_0)$-) algebra $\mathcal{B}=\bigoplus_{i\geq 0}\mathcal{B}_i$ (resp. $\mathcal{C}=\bigoplus_{i\geq 0}\mathcal{C}_i$). The morphism $q$ is dual to $\mathcal{B}\otimes_{\mathcal{B}_0}\mathcal{C}_0\to \mathcal{C}$. Consider a point $z\in\mathcal{B}\otimes_{\mathcal{B}_0}\mathcal{C}_0$ with residue field $\kappa(x)$, and consider a minimal system $f_1,...,f_r$ of homogeneous generators of the ideal $\mathcal{C}_{>0}/(\mathcal{B}_{>0})\otimes\kappa(z)$ of $\mathcal{C}/(\mathcal{B}_{>0})\otimes\kappa(z)$. Shrink $ X^0\times_{Y^0} Y^\pm$ to an affine neighborhood $z$ so that those generators lift to generators $F_1,...,F_r$ of $\mathcal{C}_{>0}/(\mathcal{B}_{>0})$. On the one hand, we have a natural morphism $p:\mathcal{B}\otimes_{\mathcal{B}_0}\mathcal{C}_0[F_1,...,F-r]\to \mathcal{C}$. Since $F_1,...,F_r$ are homogeneous and generate the ideal $\mathcal{C}/(\mathcal{B}_{>0})$, the morphism $p$ is surjective by induction on the degrees.\medskip

On the other hand, note that $f_1,...,f_r$ are linearly independent in the relative cotangent space
of $q$. As the fibers are smooth, we have $r\leq \dim(q^{-1}(z))=\dim(Spec(\mathcal{C}\otimes\kappa(z)))$. Since $q$ is flat, we have:
\begin{align}
    \mathcal{C}\otimes\kappa(z)=(\mathcal{B}\otimes_{\mathcal{B}_0}\mathcal{C}_0[F_1,...,F-r]/\ker(p))\otimes\kappa(z)=\frac{\kappa(z)[F_1,...,F_r]}{\ker(p)\otimes\kappa(z)}
\end{align}
The right hand side has dimension $r$ only if $\ker(p)\otimes\kappa(z)= 0$. Thus (after possibly shrinking $X^0\times_{Y^0} Y^\pm$
again), we have $\ker(p) = 0$ and $\mathcal{C}\simeq\mathcal{B}\otimes_{\mathcal{B}_0}\mathcal{C}_0[F_1,...,F-r]$, so $X^\pm=\mathbb{A}^r\times X^0\times_{Y^0} Y^\pm$ locally on $X^0\times_{Y^0} Y^\pm$, hence $X^\pm\to X^0\times_{Y^0} Y^\pm$ is an affine fibre bundle.$\Box$ \medskip

Denote now by $\Pi$ the set of cocharacters of $\mathbb{C}^\ast$, and for a given cocharacter $\pi\in\Pi$, by $d^+_\pi,d^0_\pi,d^-_\pi$ respectively the number of $\mathbb{C}^\ast$-contracting, resp. $\mathbb{C}^\ast$-invariants, resp. $\mathbb{C}^\ast$-repelling weights, and by $\Ind_\pi=d^+_\pi-d^-_\pi$. For a smooth $\mathbb{C}^\ast$-equivariant morphism $\phi:X\to Y$, denote by $\phi^0_\pi$ the morphism $\phi^0$ restricted to the locus of $ X^0$ where the $\mathbb{C}^\ast$-action on $X$ induces an action of cocharacter $\pi$ on $T_{X/Y}$ (in particular it is smooth of relative dimension $d^0_\pi$). It induces a decomposition into connected components, hence $\phi^0=\bigsqcup_{\pi\in\Pi}\phi^0_\pi$, and similarly $\phi^+=\bigsqcup_{\pi\in\Pi}\phi^+_\pi$.\medskip

\begin{proposition}\label{propsmoothhyploc}
    For a smooth $\mathbb{C}^\ast$-equivariant morphism $\phi:X\to Y$ of relative dimension $d$ there are natural isomorphisms in the six-functor formalism for constructible complexes:
    \begin{align}
    (p_X^\pm)_!(\eta_X^\pm)^\ast\phi^\ast[d]\simeq\bigoplus_{\pi\in\Pi}\mathbb{L}^{\pm \Ind_\pi/2}(\phi^0_\pi)^\ast[d^0_\pi](p_Y^\pm)_!(\eta_Y^\pm)^\ast
    \end{align}
    moreover, they commute with the duality $S_X$ and $S_Y$ in an obvious way.
\end{proposition}

In particular if $\phi$ is étale (resp. a local biholomorphism), $T_{X/Y}$ is trivial, hence one has a natural isomorphism:
\begin{align}
    (p_X^\pm)_!(\eta_X^\pm)^\ast\phi^\ast[d]\simeq(\phi^0)^\ast[d](p_Y^\pm)_!(\eta_Y^\pm)^\ast
    \end{align}

Proof: Consider the following commutative diagram:
\[\begin{tikzcd}\label{diagsmoothhyploc}
    X\arrow[dd,"\phi"] & X^\pm\arrow[l,"\eta_X^\pm"]\arrow[rr,"p_X^\pm"]\arrow[dd,"\phi^\pm"]\arrow[dr,"f"] && X^0\arrow[dd,"\phi^0"]\\
    && Y^\pm\times_{Y^0}X^0\arrow[dl,"\tilde{\phi}^\pm"]\arrow[ur,"\tilde{p}_X^\pm"] &\\
    Y & Y^\pm\arrow[l,"\eta_Y^\pm"]\arrow[rr,"p_Y^\pm"]& & Y^0
\end{tikzcd}\]
One has the following sequence of isomorphisms:
\begin{align}
    (p_X^\pm)_!(\eta_X^\pm)^\ast\phi^\ast[d]\simeq& (p_X^\pm)_!(\phi^\pm)^\ast[d](\eta_Y^\pm)^\ast\nn\\
    \simeq &(\tilde{p}_X^\pm)_!f_!f^\ast(\tilde{\phi}^\pm)^\ast[d](\eta_Y^\pm)^\ast
\end{align}
where we have used the commutativity of the left square in the first line, and the commutativity of the two  small triangles in the second line. The right square of the above diagram admits an decomposition as a disjoint union of those diagrams:
\[\begin{tikzcd}
    X^\pm_\pi\arrow[rr,"p_{X,\pi}^\pm"]\arrow[dd,"\phi^\pm_\pi"]\arrow[dr,"f_\pi"] && X^0_\pi\arrow[dd,"\phi^0_\pi"]\\
    & Y^\pm\times_{Y^0}X^0_\pi\arrow[dl,"\tilde{\phi}^\pm_\pi"]\arrow[ur,"\tilde{p}_{X,\pi}^\pm"] &\\
    Y^\pm\arrow[rr,"p_Y^\pm"]& & Y^0
\end{tikzcd}\]
hence one obtains:
\begin{align}
    (\tilde{p}_{X,\pi}^\pm)_!(f_\pi)_!(f_\pi)^\ast(\tilde{\phi}^\pm_\pi)^\ast[d](\eta_Y^\pm)^\ast\simeq&(\tilde{p}_{X,\pi}^\pm)_![-2d^\pm_\pi](\tilde{\phi}^\pm_\pi)^\ast[d^+_\pi+d^0_\pi+d^-_\pi](\eta_Y^\pm)^\ast\nn\\
    \simeq&\mathbb{L}^{\pm \Ind_\pi/2}\phi^0_\pi[d^0_\pi](p^\pm_Y)_!(\eta_Y^\pm)^\ast
\end{align}
where we have used the fact that $f_\pi$ is an affine fibration with $d^\pm_\pi$-dimensional fiber in the first line and the base change in the down right Cartesian triangle in the second line. Summing on the connected components decomposition, one obtains a natural isomorphism in the six-functor formalism:
\begin{align}
    (p_X^\pm)_!(\eta_X^\pm)^\ast\phi^\ast[d]\simeq \bigoplus_{\pi\in\Pi}\mathbb{L}^{\pm \Ind_\pi/2}(\phi^0)^\ast[d^0_\pi](p_Y^\pm)_!(\eta_Y^\pm)^\ast
\end{align}
$\Box$\medskip

Now considering $U$ a smooth $S$-algebraic space with $\mathbb{C}^\ast$-action, one can consider the smooth $\mathbb{C}^\ast$-equivariant morphism $\phi:U\to S$, where $S$ has the trivial $\mathbb{C}^\ast$-action. The hyperbolic localization functor is then trivial on $S$, and:
\begin{align}
   \mathcal{IC}_U&=\phi^\ast[d]\mathcal{IC}_S=\phi^![-d]\mathcal{IC}_S
\end{align}
in particular $\mathcal{IC}_S$ is $\mathbb{C}^\ast$-equivariant, hence $\mathcal{IC}_U$ is also $\mathbb{C}^\ast$-equivariant, and self-Verdier dual, hence the natural arrow $\mathbb{D}(p_U^-)_!(\eta_U^-)^\ast\mathcal{IC}_U\to(p_X^+)_!(\eta_X^+)^\ast\mathcal{IC}_U$ is an isomorphism. Moreover, there is then an open and closed decomposition $U^0=\bigsqcup_{\pi\in\Pi}U^0_\pi$, and:
\begin{align}
   \mathcal{IC}_{U^0_\pi}&=(\phi^0_\pi)^\ast[d^0_\pi]\mathcal{IC}_S=(\phi^0_\pi)^![-d^0_\pi]\mathcal{IC}_S
\end{align}
hence an absolute version of the last proposition gives:

\begin{corollary}\label{classbibi}
    For a smooth $\mathbb{C}^\ast$-equivariant $S$-algebraic (resp. analytic) space $U$ there are natural isomorphisms in the six-functor formalism such that the following diagram commutes:
    \[\begin{tikzcd}
        \mathbb{D}(p_U^-)_!(\eta_U^-)^\ast\mathcal{IC}_U\arrow[dd,"S_U"]\arrow[dr,"\simeq"]&\\ &\bigoplus_{\pi\in\Pi}\mathbb{L}^{\Ind_\pi/2}\mathcal{IC}_{U^0_\pi}\\
        (p_U^+)_!(\eta_U^+)^\ast\mathcal{IC}_U\arrow[ur,"\simeq"] &
    \end{tikzcd}\]
\end{corollary}

\section{Critical Białynicki-Birula decomposition}

In this section and in the rest of the papers, all the algebraic spaces will be assumed to be quasi-separated and locally of finite type over $\mathbb{C}$: in particular, all the $\mathbb{C}^\ast$-actions will then be étale locally linearizable, the proper pullback and proper pushout of any functor, and the hyperbolic localization diagram, will be defined.

\subsection{The functor of vanishing cycles}

Consider an algebraic space $U$ with a regular function $f:U\to\mathbb{C}$. We give a definition of the functor of vanishing cycles $\phi_f$, which is equivalent to those exposed in definition 2.10 of \cite{Joycesymstab}. Consider the following commutative diagram of complex analytic spaces, where $U^{an}$ is the analytification of $U$, and the square are Cartesian:
\[\begin{tikzcd}
    U^{an}_0\arrow[r,"\iota"]\arrow[d,"f"] & U^{an}\arrow[d,"f"] & \tilde{U}^{an}_\theta\arrow[l,swap,"j_\theta"]\arrow[d,"f"]\\
    \{0\}\arrow[r] & \mathbb{C} & e^{i\theta}\{v\in\mathbb{C}|\Re(v)\leq 0\}\arrow[l]
\end{tikzcd}\]
We define the functors:
\begin{align}
    \phi_f^\theta:=\iota^\ast(j_\theta)_!(j_\theta)^!:D(X)\to D(X_0)
\end{align}
There are natural isomorphisms $\tilde{T}^\theta_f:\phi^0_f\to\phi^\theta_f$ induced by rotating $\mathbb{C}$ counterclockwise by an angle $\theta$ around the origin. One defines the vanishing cycles functor $\phi_f:=\phi^0_f$, and the natural isomorphism $\tilde{T}^{2\pi}_f:\phi_f\to\phi_f$ is the monodromy of vanishing cycle. Because $\phi_f$ is constructed using maps of analytic spaces and not maps of $\mathbb{C}$-schemes, it is not direct that they maps constructible complex of $U$ to constructible complex of $U_0$, but it is the case, as explained in definition 2.10 of \cite{Joycesymstab}, \ie it defines a functor $\phi_f:D^b_c(U)\to D^b_c(U_0)$, and it even defines functors $\phi_f:\Perv(U)\to \Perv(U_0)$.\medskip

Massey defines in \cite{Massey2016NaturalCO} a natural isomorphism giving the commutation of the vanishing cycles functor (considered as a functor $\phi_f:D^b_c(U)\to D^b_c(U_0)$) with the Verdier duality $\mathbb{D}\phi_f\simeq \phi_f\mathbb{D}$. Consider the Cartesian diagram of closed subspaces:
\[\begin{tikzcd}
    U^{an}_0\arrow[dr,"q"] & &\\
    &\tilde{U}^{an}_0\cap\tilde{U}^{an}_\pi\arrow[r,"\hat{j}_\pi"]\arrow[d,"\hat{j}_0"]\arrow[dr,"m"] & \tilde{U}^{an}_0\arrow[d,"j_0"]\\
    &\tilde{U}^{an}_\pi\arrow[r,"j_\pi"] & U
\end{tikzcd}\]
Then Massey defines the following sequence of isomorphisms:
\begin{align}
    \mathbb{D}\phi_f\simeq&\mathbb{D}q^\ast m^\ast (j_0)_!(j_0)^!\nn\\
    \simeq&\mathbb{D}q^!m^\ast (j_0)_!(j_0)^!\nn\\
    \simeq&\mathbb{D}q^!(\hat{j}_0)^!(j_\pi)^\ast\nn\\
    \simeq&q^\ast(\hat{j}_0)^\ast(j_\pi)^!\mathbb{D}\nn\\
    \simeq&q^\ast m^\ast(j_\pi)_!(j_\pi)^!\mathbb{D}\nn\\
    \simeq&\phi^\pi_f\mathbb{D}
\end{align}
where in the first and last line we have used the fact that $\iota=m\circ q$, in the second line \cite[Lem 2.1]{Massey2016NaturalCO} which proves that the natural morphism $q^!m^\ast (j_0)_!(j_0)^!\to q^\ast m^\ast (j_0)_!(j_0)^!$ is an isomorphism and in the third and fifth line the isomorphisms built in \cite[Lem 2.2]{Massey2016NaturalCO} using the six-functor formalism in the Cartesian square of closed subspace of $X$, Using the fact that $U^{an}=\tilde{U}^{an}_0\cup\tilde{U}^{an}_\pi$. The Massey isomorphism is then defined by the composition:
\[\begin{tikzcd}
    \mathbb{D}\phi_f\arrow[r,"\simeq"] & \phi^\pi_f\mathbb{D}\arrow[r,"(\tilde{T}_f^\pi)^{-1}"] & \phi_f\mathbb{D} 
\end{tikzcd}\]
and it commutes with the monodromy operator $\tilde{T}^{2\pi}_f$.\medskip

For two algebraic spaces $X,Y$ with projections $\pi_1:X\times Y\to X,\pi_2:X\times Y\to Y$, for $F\in D^b_c(X),G\in D^b_c(Y)$, one denotes:
\begin{align}
    F\boxtimes G:=\pi_1^\ast(F)\otimes\pi_2^\ast G
\end{align}
For $f;X\to\mathbb{C},g:Y\to\mathbb{C}$, denote by $f\boxplus g:=f\circ\pi_1+g\circ\pi_2:X\times Y\to\mathbb{C}$. For $j:Y\to X$ a closed embedding, denote $R\Gamma_Y:=j_! j^!:D^b_c(X)\to D^b_c(X)$.\medskip

Consider two algebraic spaces $U,V$ with a regular function $f:U\to\mathbb{C},g:V\to\mathbb{C}$, denote by $k:U_0\times V_0\to(U\times V)_0$ the embedding into the zero locus of $f\boxplus g$. Massey defines in \cite{MasThomSeb} a natural Thom-Sebastiani isomorphism:
\begin{align}
    \phi_f\boxtimes\phi_g\simeq k^\ast\phi_{f\boxplus g}
\end{align}
denote by $\iota:U_0\to U,\jmath:V_0\to V,q:(U\times V)_0\to U\times V$ the closed embeddings of the zero set of $f,g,f\boxplus g$. One has the closed embedding $\tilde{U}_0\times \tilde{V}_0\to (\tilde{U\times V})_0$, and \cite[lem 1.2]{MasThomSeb} shows that the natural arrow of the six-functor formalism $R\Gamma_{\tilde{U}_0\times \tilde{V}_0}\to R\Gamma_{(\tilde{U\times V})_0}$ is an isomorphism. The Thom-Sebastiani isomorphism is then defined by the composition of natural isomorphisms of the six-functor formalism:
\begin{align}
    \phi_f\boxtimes\phi_g=&\iota^\ast R\Gamma_{\tilde{U}_0}\boxtimes \jmath^\ast R\Gamma_{\tilde{V}_0}\nn\\
    \simeq&(\iota\times\jmath)^\ast R\Gamma_{\tilde{U}_0\times \tilde{V}_0}\nn\\
    \overset{\simeq}{\to}& k^\ast q^\ast R\Gamma_{\tilde{U\times V}_0}\nn\\
    =&k^\ast\phi_{f\boxplus g}
\end{align}
it is shown in \cite{MasThomSeb} that this isomorphism commutes with monodromy and it can be shown that it commutes with duality.\medskip

There is a theory of vanishing cycles functor for mixed Hodge modules, which projects under $rat$ to the above theory for perverse sheaves, as exposed in \cite[sec 2.10]{Joycesymstab}. There is a vanishing cycle functor:
\begin{align}
    \phi^H_f:\MHM(U)\to \MMHM(U_0)
\end{align}
It sends a mixed Hodge module on $U$ to a monodromic mixed Hodge module on $U_0$, \ie a mixed Hodge module with commuting actions of a unipotent operator $T_s$ and a nilpotent operator $N$, giving respectively the semisimple part and the logarithm of the unipotent part of the monodromy operator. There is also a self duality isomorphism and a Thom-Sebastiani isomorphism, and it is checked in \cite[appendix A]{Joycesymstab} that $rat$ sends them to the corresponding isomorphisms of perverse sheaves. When a mixed Hodge module $M\in \MHM(U)$ is polarized, \ie is provided with an isomorphism $\sigma:\mathbb{D}M\simeq M$, the self-duality isomorphism of $\phi^H_f$ provides them $\phi^H_fM$ with a strong polarization, \ie an isomorphism with its dual commuting with monodromy.

\subsection{Functoriality of the vanishing cycles}

The vanishing cycles functor has some functoriality properties. Consider a map of algebraic spaces $\Phi:U\to V$, and a regular functions $f:V\to\mathbb{C}$, $f=g\circ\Phi:U\to\mathbb{C}$, and denote by $\Phi_0:U_0\to V_0$ the induced map on the zero locus.

\begin{proposition}\label{propfunctvan}
    There are natural morphisms built using the six-functor formalism for constructible complexes or complexes of mixed Hodge modules:
    \begin{align}
         (\Phi_0)^\ast\phi_g^\theta&\to\phi_f^\theta\Phi^\ast\nn\\
         (\Phi_0)_!\phi_f^\theta&\to\phi_g^\theta\Phi_!
    \end{align}
    They are compatible with monodromy and Thom-Sebastiani isomorphism, and are compatible with composition and base change of $\Phi$.
\end{proposition}

Proof: We consider then the commutative diagram, where all the squares are Cartesian:
\[\begin{tikzcd}
    U^{an}_0\arrow[r,"\iota_U"]\arrow[d,"\Phi_0"] & U^{an}\arrow[d,"\Phi"] & \tilde{U}^{an}_\theta\arrow[l,"j_ {\theta,U}"]\arrow[d,"\tilde{\Phi}_\theta"]\\
    V^{an}_0\arrow[r,"\iota_V"]\arrow[d,"g"] & V^{an}\arrow[d,"g"] & \tilde{V}^{an}_\theta\arrow[l,"j_{\theta,V}"]\arrow[d,"g"]\\
    \{0\}\arrow[r] & \mathbb{C} & e^{i\theta}\{v\in\mathbb{C}|\Re(v)\leq 0\}\arrow[l]
\end{tikzcd}\]
In the following of the proof, the superscript $"an"$ will be implicit for readability. Then one has the following sequence of natural morphisms of the six-functor formalism:
\begin{align}
    (\Phi_0)^\ast\phi_g^\theta=&\Phi_0^\ast(\iota_V)^\ast(j_{\theta,V})_!(j_{\theta,V})^!\nn\\
    \simeq&(\iota_U)^\ast\Phi^\ast(j_{\theta,V})_!(j_{\theta,V})^!\nn\\
    \simeq&(\iota_U)^\ast(j_{\theta,U})_!\tilde{\Phi}^\ast(j_{\theta,V})^!\nn\\
    \to&(\iota_U)^\ast(j_{\theta,U})_!(j_{\theta,U})^!\Phi^\ast\nn\\
    =&\phi_f^\theta\Phi^\ast
\end{align}
where the first and last lines are the definitions, the second line follows from the commutativity of the left square, the third from the Cartesianity of the right square, and then fourth from the commutativity of the right square.

\begin{align}
    (\Phi_0)_!\phi_f^\theta=&(\Phi_0)_!(\iota_U)^\ast(j_{\theta,U})_!(j_{\theta,U})^!\nn\\
    \simeq&(\iota_V)^\ast\Phi_!(j_{\theta,U})_!(j_{\theta,U})^!\nn\\
    \simeq&(\iota_V)^\ast(j_{\theta,VU})^!\nn\\
    \to&(\iota_V)^\ast(j_{\theta,V})_!(j_{\theta,V})^!\Phi_!\nn\\
    =&\phi_g^\theta\Phi_!
\end{align}
where the first and last lines are the definitions, the second line follows from the Cartesianity of the left square, the third from the commutativity of the right square, and then fourth from the fact that the right square is Cartesian.\medskip

Because $f=g\circ\Phi$ these morphisms are compatible with the isomorphisms $\tilde{T}^\theta$, namely the following squares are commutative:
\[\begin{tikzcd}\label{diagfunctmon}
    (\Phi_0)^\ast\phi_g\arrow[r]\arrow[d,"(\Phi_0)^\ast\tilde{T}^\theta_g"] & \phi_f\Phi^\ast\arrow[d,"\tilde{T}^\theta_f\Phi^\ast"] & (\Phi_0)_!\phi_f\arrow[r]\arrow[d,"(\Phi_0)_!\tilde{T}^\theta_f"] & \phi_g\Phi^\ast\arrow[d,"\tilde{T}^\theta_g\Phi_!"]\\
    (\Phi_0)^\ast\phi_g^\theta\arrow[r] & \phi_f^\theta\Phi^\ast & (\Phi_0)_!\phi_f^\theta\arrow[r] & \phi_g^\theta\Phi_!
\end{tikzcd}\]
in particular these morphisms are compatible with the monodromy.\medskip

Consider now $U,U',V,V'$ algebraic spaces with regular functions to $\mathbb{C}$ $f,f',g,g'$ and maps $\Phi:U\to U'$, $\Psi:V\to V'$. Considering the Cartesian diagram, where the horizontal maps are closed embeddings:
\[\begin{tikzcd}
    U\times V\arrow[d,"\Phi\times\Psi"] & (\tilde{U\times V})_0\arrow[d,"\Phi\times\Psi"]\arrow[l] & \tilde{U}_0\times\tilde{V}_0\arrow[d,"\Phi\times\Psi"]\arrow[l]\\
    U'\times V' & (\tilde{U'\times V'})_0\arrow[l] & \tilde{U'}_0\times\tilde{V'}_0\arrow[l]\\
\end{tikzcd}\]
One has then commutative squares of morphisms from the six-functor formalism:
\[\begin{tikzcd}[column sep= 2ex]
    (\Phi\times\Psi)^\ast R\Gamma_{\tilde{U'}_0\times\tilde{V'}_0}\arrow[r]\arrow[d,"\simeq"] & R\Gamma_{\tilde{U}_0\times\tilde{V}_0}(\Phi\times\Psi)^\ast\arrow[d,"\simeq"] & (\Phi\times\Psi)_! R\Gamma_{\tilde{U}_0\times\tilde{V}_0}\arrow[r]\arrow[d,"\simeq"] & R\Gamma_{\tilde{U'}_0\times\tilde{V'}_0}(\Phi\times\Psi)_!\arrow[d,"\simeq"] 
    \\
    (\Phi\times\Psi)^\ast R\Gamma_{(\tilde{U'\times V'})_0}\arrow[r] & R\Gamma_{(\tilde{U\times V})_0}(\Phi\times\Psi)^\ast & (\Phi\times\Psi)_! R\Gamma_{(\tilde{U\times V})_0}\arrow[r] & R\Gamma_{(\tilde{U'\times V'})_0}(\Phi\times\Psi)_!
\end{tikzcd}\]
Consider now the Cartesian diagram:
\[\begin{tikzcd}
    U_0\times V_0\arrow[r,"k"]\arrow[d,"\Phi_0\times\Psi_0"] & (U\times V)_0\arrow[r,"q"]\arrow[d,"(\Phi\times\Psi)_0"] & U\times V\arrow[d,"\Phi\times\Psi"]\\
    U'_0\times V'_0\arrow[r,"k"] & (U'\times V')_0\arrow[r,"q"] & U'\times V'
\end{tikzcd}\]
One has then commutative squares of morphisms from the six-functor formalism:
\[\begin{tikzcd}[column sep= 2ex]
    (\Phi_0\times\Psi_0)^\ast(\iota'\times\jmath')^\ast\arrow[r]\arrow[d,"\simeq"] & (\iota\times\jmath)^\ast(\Phi\times\Psi)^\ast\arrow[d,"\simeq"] & (\Phi_0\times\Psi_0)_!(\iota\times\jmath)^\ast\arrow[r]\arrow[d,"\simeq"] & (\iota'\times\jmath')^\ast(\Phi\times\Psi)_!\arrow[d,"\simeq"]\\
    k^\ast((\Phi\times\Psi)_0)^\ast (q')^\ast\arrow[r] & (\iota\times\jmath)^\ast(\Phi\times\Psi)^\ast & k^\ast((\Phi\times\Psi)_0)_! q^\ast\arrow[r] & (\iota'\times\jmath')^\ast(\Phi\times\Psi)_!
\end{tikzcd}\]
And one obtains then finally commutative square of morphisms:
\[\begin{tikzcd}[column sep= 2ex]
     (\Phi_0\times\Psi_0)^\ast(\iota'\times\jmath')^\ast R\Gamma_{\tilde{U'}_0\times\tilde{V'}_0}\arrow[r,"\simeq"]\arrow[d,"\simeq"] & (\iota\times\jmath)^\ast(\Phi\times\Psi)^\ast R\Gamma_{\tilde{U'}_0\times\tilde{V'}_0}\arrow[r]\arrow[d,"\simeq"] & (\iota\times\jmath)^\ast R\Gamma_{\tilde{U}_0\times\tilde{V}_0}(\Phi\times\Psi)^\ast\arrow[d,"\simeq"]\\
     k^\ast(\Phi\times\Psi)_0^\ast (q')^\ast R\Gamma_{\tilde{U'}_0\times\tilde{V'}_0}\arrow[r,"\simeq"]\arrow[d,"\simeq"] & (\iota\times\jmath)^\ast(\Phi\times\Psi)^\ast R\Gamma_{\tilde{U'}_0\times\tilde{V'}_0} \arrow[r]\arrow[d,"\simeq"] & (\iota\times\jmath)^\ast R\Gamma_{\tilde{U}_0\times\tilde{V}_0}(\Phi\times\Psi)^\ast\arrow[d,"\simeq"]\\
     k^\ast(\Phi\times\Psi)_0^\ast (q')^\ast R\Gamma_{\tilde{U'}_0\times\tilde{V'}_0}\arrow[r,"\simeq"] & (\iota\times\jmath)^\ast(\Phi\times\Psi)^\ast R\Gamma_{\tilde{U'}_0\times\tilde{V'}_0} \arrow[r] & (\iota\times\jmath)^\ast R\Gamma_{\tilde{U}_0\times\tilde{V}_0}(\Phi\times\Psi)^\ast
\end{tikzcd}\]
And similarly:
\[\begin{tikzcd}[column sep= 2ex]
     (\Phi_0\times\Psi_0)_!(\iota\times\jmath)^\ast R\Gamma_{\tilde{U}_0\times\tilde{V}_0}\arrow[r,"\simeq"]\arrow[d,"\simeq"] & (\iota'\times\jmath')^\ast(\Phi\times\Psi)_! R\Gamma_{\tilde{U}_0\times\tilde{V}_0}\arrow[r]\arrow[d,"\simeq"] & (\iota'\times\jmath')^\ast R\Gamma_{\tilde{U'}_0\times\tilde{V'}_0}(\Phi\times\Psi)_!\arrow[d,"\simeq"]\\
     (k')^\ast((\Phi\times\Psi)_0)_! q^\ast R\Gamma_{\tilde{U}_0\times\tilde{V}_0}\arrow[r,"\simeq"]\arrow[d,"\simeq"] & (\iota'\times\jmath')^\ast(\Phi\times\Psi)_! R\Gamma_{\tilde{U}_0\times\tilde{V}_0} \arrow[r]\arrow[d,"\simeq"] & (\iota'\times\jmath')^\ast R\Gamma_{\tilde{U'}_0\times\tilde{V'}_0}(\Phi\times\Psi)_!\arrow[d,"\simeq"]\\
     (k')^\ast((\Phi\times\Psi)_0)_! q^\ast R\Gamma_{\tilde{U}_0\times\tilde{V}_0}\arrow[r,"\simeq"] & (\iota'\times\jmath')^\ast(\Phi\times\Psi)_! R\Gamma_{\tilde{U}_0\times\tilde{V}_0} \arrow[r] & (\iota'\times\jmath')^\ast R\Gamma_{\tilde{U'}_0\times\tilde{V'}_0}(\Phi\times\Psi)_!
\end{tikzcd}\]
Hence the following diagrams are commutative:
\[\begin{tikzcd}[column sep= 1ex]\label{diagcomthomseb}
    (\Phi_0\times\Psi_0)^\ast(\phi_{f'}\boxtimes\phi_{g'})\arrow[r]\arrow[d,"\simeq"] &  (\phi_f\boxtimes\phi_g)(\phi\times\Psi)^\ast\arrow[d,"\simeq"] & (\Phi_0\times\Psi_0)_!(\phi_f\boxtimes\phi_g)\arrow[r]\arrow[d,"\simeq"] &  k^\ast(\phi_{f'}\boxtimes\phi_{g'})(\phi\times\Psi)_!\arrow[d,"\simeq"]\\
    k^\ast((\Phi\times\Psi)_0)^\ast\phi_{f'\boxplus g'}\arrow[r] & \phi_{f\boxplus g}(\Phi\times\Psi)^\ast & (k')^\ast((\Phi\times\Psi)_0)_!\phi_{f\boxplus g}\arrow[r] & k^\ast\phi_{f'\boxplus g'}(\Phi\times\Psi)_!
\end{tikzcd}\]
\ie the morphisms of the proposition are compatible with the Thom-Sebastiani isomorphism.\medskip

Consider now $U,V,W$ algebraic spaces, with three regular maps to $\mathbb{C}$, and maps $\Phi:U\to V,\Psi:V\to W$ such that $f=g\circ\Phi, g=h\circ\Psi$. Considering the naturality of the six-functor formalism in the Cartesian diagram:
\[\begin{tikzcd}
    U_0\arrow[r,"\iota_U"]\arrow[d,"\Phi_0"] & U\arrow[d,"\Phi"] & \tilde{U}_\theta\arrow[l,"j_ {\theta,U}"]\arrow[d,"\tilde{\Phi}_\theta"]\\
    V_0\arrow[r,"\iota_V"]\arrow[d,"\Psi_0"] & V\arrow[d,"\Psi"] & \tilde{V}_\theta\arrow[l,"j_ {\theta,V}"]\arrow[d,"\tilde{\Psi}_\theta"]\\
    W_0\arrow[r,"\iota_W"]\arrow[d,"h"] & W\arrow[d,"h"] & \tilde{W}_\theta\arrow[l,"j_{\theta,W}"]\arrow[d,"h"]\\
    \{0\}\arrow[r] & \mathbb{C} & e^{i\theta}\{v\in\mathbb{C}|\Re(v)\leq 0\}\arrow[l]
\end{tikzcd}\]
one obtains directly that the following diagrams of morphisms are commutative:
\[\begin{tikzcd}
    & (\Phi_0)^\ast\phi_g\Psi^\ast\arrow[dr] & \\
    ((\Psi\circ\Phi)_0)^\ast\phi_h=(\Phi_0)^\ast(\Psi_0)^\ast\phi_h\arrow[ur]\arrow[rr] && \phi_f\Phi^\ast\Psi^\ast=\phi_f(\Psi\circ\Phi)^\ast\\
    & (\Psi_0)_!\phi_g\Phi_!\arrow[dr] & \\
    ((\Psi\circ\Phi)_0)_!\phi_f=(\Psi_0)_!(\Phi_0)_!\phi_f\arrow[ur]\arrow[rr] && \phi_f\Psi_!\Phi_!=\phi_f(\Psi\circ\Phi)_!
\end{tikzcd}\]
hence the morphisms of the proposition are compatible with the composition.\medskip

Consider now a Cartesian square of morphisms between algebraic spaces:
\[\begin{tikzcd}
    U\times_W V\arrow[r,"\Phi'"]\arrow[d,"\Psi'"] & V\arrow[d,"\Psi"]\\
    U\arrow[r,"\Phi"] & W
\end{tikzcd}\]
and a regular function $f:W\to\mathbb{C}$. Consider now the following commutative diagram, where all the squares are Cartesian:
\[\begin{tikzcd}[column sep= 1ex]
    & (U\times_W V)_0\arrow[rrr]\arrow[dl]\arrow[ddr] &&& U\times_W V\arrow[dl]\arrow[ddr] &&& (\tilde{U\times_W V})_0\arrow[lll]\arrow[dl]\arrow[ddr]&\\
    U_0\arrow[rrr]\arrow[ddr] &&& U\arrow[ddr] &&& \tilde{U}_0\arrow[lll]\arrow[ddr]&&\\
    &&V_0\arrow[rrr]\arrow[dl] &&& V\arrow[dl] &&& \tilde{V}_0\arrow[lll]\arrow[dl]\\
    &W_0\arrow[rrr] &&& W &&& \tilde{W}_0\arrow[lll]&\\
\end{tikzcd}\]
Hence from the naturality of the isomorphisms expressing the functoriality and base change of the six-functor formalism, one obtains that the following diagram of morphisms is commutative:
\[\begin{tikzcd}
    & (\Phi'_0)_!(\Psi'_0)^\ast\phi_{f\circ\Phi}\arrow[r] & (\Phi'_0)_!\phi_{f\circ\Phi\circ\Psi'}(\Psi')^\ast\arrow[d]\\
    (\Psi_0)^\ast(\Phi_0)_!\phi_{f\circ\Phi}\arrow[d]\arrow[ur,"\simeq"] && \phi_{f\circ\Psi}(\Phi')_!(\Psi')^\ast\\
    (\Psi_0)^\ast\phi_{f}\Phi_!\arrow[r] & \phi_{f\circ\Psi}\Psi^\ast\Phi_!\arrow[ur,"\simeq"]
\end{tikzcd}\]
where the diagonal arrows are the isomorphisms from the base change and the vertical and horizontal arrows are the morphisms from the proposition. Then the morphisms of the proposition are compatible with base change.\medskip

The functoriality of functor of vanishing cycles for mixed Hodge modules shows that the morphisms of the proposition are defined at the level of mixed Hodge modules. Because of the compatibility results of \cite[Appendix A]{Joycesymstab}, the functor $rat$ sends the various squares expressing the compatibility conditions in this proposition to the corresponding squares at the level of perverse sheaves, which are commutative, hence they commute at the level of mixed Hodge modules.
$\Box$\medskip

Now, for $\Phi$ smooth of relative dimension $d$ (resp. $\Psi$ proper), observing that $\mathbb{D}\Phi^\ast[d]\mathbb{D}=\Phi^\ast[d]$ (resp. $\mathbb{D}\Psi_!\mathbb{D}=\Psi_!$ and then also $\mathbb{D}(\Phi_0)^\ast[d]\mathbb{D}=(\Phi_0)^\ast[d]$ because $\Phi_0$ is smooth of codimension $d$ (resp. $\mathbb{D}(\Psi_0)_!\mathbb{D}=(\Psi_0)_!$ because $\Psi_0$ is proper), consider the following diagrams:
\[\begin{tikzcd}
    (\Phi_0)^\ast[d]\phi_g\arrow[r]\arrow[d,"\simeq"] & \phi_f\Phi^\ast[d] \arrow[d,"\simeq"] & (\Psi_0)_!\phi_f\arrow[r]\arrow[d,"\simeq"] & \phi_g\Psi_! \arrow[d,"\simeq"]\\
    \mathbb{D}(\Phi_0)^\ast[d]\phi_g\mathbb{D} & \mathbb{D}\phi_f\Phi^\ast[d]\mathbb{D}\arrow[l] & \mathbb{D}(\Psi_0)_!\phi_f\mathbb{D} & \mathbb{D}\phi_g\Psi_!\mathbb{D}\arrow[l]
\end{tikzcd}\]
where the vertical arrows comes from the self-duality of the vanishing cycles functor. It is not hard to show that these squares are commutative, and we find again the standard fact that for $\Phi$ smooth of relative dimension $d$, and $\Psi$ proper, there are natural isomorphisms:
\begin{align}
    (\Phi_0)^\ast[d]\phi_g\simeq &\phi_f\Phi^\ast[d]\nn\\
    (\Psi_0)_!\phi_f\simeq &\phi_g\Psi_!
\end{align}
commuting with the monodromy and the duality of the vanishing cycles.\medskip

Consider now a $\mathbb{C}^\ast$ action on $U$ such that $f$ is $\mathbb{C}^\ast$ invariant, \ie $U$ is can be considered as an algebraic space over $\mathbb{A}^1_\mathbb{C}$. In particular $U$ is locally of finite type and quasi-separated over a quasi-separated base, hence the $\mathbb{C}^\ast$-action on it is étale-locally linearizable, and then the $\mathbb{C}^\ast$-action on $U^{an}$ (and also on any base change of $U^{an}$) is analytic-locally linearizable over $\mathbb{C}$, hence the hyperbolic localization diagram of $U^{an}$ and all its base change exists as diagrams of analytic spaces. Consider a $\mathbb{C}^\ast$-equivariant constructible complex $A$. The functor $\phi_f$ is defined using the morphisms of the six-functor formalism defined $\mathbb{C}^\ast$-equivariant functions, hence $\phi_f(A)$ is $\mathbb{C}^\ast$-equivariant. We will use the idea of \cite{Ric16} to prove the commutation of the hyperbolic localization with the vanishing cycle:\medskip

\begin{proposition}\label{comhyplocvan}
    For $A\in D^b_c(U)^{T-mon}$, the natural morphism $(p^\pm_{U_0})_!(\eta^\pm_{U_0})^\ast\phi_f\simeq \phi_{f^0}(p^\pm_U)_!(\eta^\pm_U)^\ast(A)$ is an isomorphism, compatible with monodromy and duality, in the sense that the following diagram is commutative:
    \[\begin{tikzcd}\label{diagcomvandual}
    (p^+_{U_0})_!(\eta^+_{U_0})^\ast\phi_f\mathbb{D}(A)\arrow[r,"\simeq"] & \phi_{f^0}(p^+_U)_!(\eta^+_U)^\ast\mathbb{D}(A) \\
    \mathbb{D}(p^-_{U_0})_!(\eta^-_{U_0})^\ast\phi_f\mathbb{D}(A)\arrow[u,"\simeq"] & \mathbb{D}\phi_{f^0}(p^-_U)_!(\eta^-_U)^\ast\mathbb{D}(A)\arrow[l,"\simeq"]\arrow[u,"\simeq"]
\end{tikzcd}\]
where the vertical arrows are given by the isomorphisms of the vanishing cycles with its dual and the isomorphism of proposition \ref{hyplocdual}, and commutes with Thom-Sebastiani isomorphism and pullback by smooth maps. 
\end{proposition}

Proof: It follows directly from proposition \ref{propfunctvan} that the isomorphism of the proposition is compatible with monodromy. To show that the natural morphism is an isomorphism, we need only to show that the diagram of the proposition is commutative. In fact, the following diagram of $\mathbb{C}^\ast$-equivariant spaces:
\[\begin{tikzcd}
    U^{an}_0\arrow[dr,"q"] & &\\
    &\tilde{U}^{an}_0\cap\tilde{U}^{an}_\pi\arrow[r,"\hat{j}_\pi"]\arrow[d,"\hat{j}_0"]\arrow[dr,"m"] & \tilde{U}^{an}_0\arrow[d,"j_0"]\\
    &\tilde{U}^{an}_\pi\arrow[r,"j_\pi"] & U^{an}
\end{tikzcd}\]
it is obtain by base change from the space $U^{an}\to 0$, hence the hyperbolic localization diagram of all the spaces in this diagram are also obtained by base change from the hyperbolic localization diagram. Hence all the squares are Cartesian in the diagram obtained by superposing the hyperbolic localization diagrams of those spaces. hence the natural isomorphisms of the six-functor formalism gives that the following square of isomorphisms is commutative:
\[\begin{tikzcd}
    (p^+_{U_0})_!(\eta^+_{U_0})^\ast\phi_f^\pi\mathbb{D}\arrow[r] & \phi_{f^0}^\pi(p^+_U)_!(\eta^+_U)^\ast\mathbb{D} \\
    \mathbb{D}(p^-_{U_0})_!(\eta^-_{U_0})^\ast\phi_f\mathbb{D}(A)\arrow[u] & \mathbb{D}\phi_{f^0}(p^-_U)_!(\eta^-_U)^\ast\mathbb{D}(A)\arrow[l]\arrow[u]
\end{tikzcd}\]
where the vertical arrow are given by Richarz morphism $S_U$ and Massey's isomorphism $\mathbb{D}\phi_f\simeq\phi_f^\pi\mathbb{D}$. Because the isomorphism of the proposition is compatible with the monodromy, hence the following square is commutative:
\[\begin{tikzcd}
    (p^+_{U_0})_!(\eta^+_{U_0})^\ast\phi_f\mathbb{D}\arrow[r]\arrow[d,"(p^+_{U_0})_!(\eta^+_{U_0})^\ast\tilde{T}^\pi_f"] & \phi_{f^0}(p^+_U)_!(\eta^+_U)^\ast\mathbb{D}\arrow[d,"\tilde{T}^\pi_{f^0}"] \\
    (p^+_{U_0})_!(\eta^+_{U_0})^\ast\phi_f^\pi\mathbb{D}\arrow[r] & \phi_f^\pi(p^+_U)_!(\eta^+_U)^\ast\mathbb{D}
\end{tikzcd}\]
and then gluing vertically the two last squares we obtain that the square of the proposition is commutative.\medskip

Consider now $U,V$ two algebraic  spaces with $\mathbb{C}^\ast$ action, $\mathbb{C}^\ast$-invariants regular functions $f:U\to\mathbb{C},g:V\to\mathbb{C}$, and a smooth map $\Phi:U\to V$ of relative dimension $d$. Applying the compatibility of the morphisms of proposition \ref{propfunctvan} with composition and base change to the morphisms in the following commutative diagram \ref{diagsmoothhyploc}, one obtains that the following diagram is commutative:
\[\begin{tikzcd}[scale cd=0.8]\label{diagcomsmoothvanhyp}
    &\bigoplus_{\pi\in\Pi}\mathbb{L}^{\pm \Ind_\pi/2}(\Phi^0_{0,\pi})^\ast[d^0_\pi](p^\pm_{V_0})_!(\eta^\pm_{V_0})^\ast\phi_g\arrow[r] &  \bigoplus_{\pi\in\Pi}\mathbb{L}^{\pm \Ind_\pi/2}(\Phi^0_{0,\pi})^\ast[d^0_\pi]\phi_{g^0}(p^\pm_V)_!(\eta^\pm_V)^\ast\arrow[d,"\simeq"] \\
    (p^\pm_{U_0})_!(\eta^\pm_{U_0})^\ast(\Phi_0)^\ast[d]\phi_g\arrow[ur,"\simeq"]\arrow[d,"\simeq"] && \phi_{f^0}\bigoplus_{\pi\in\Pi}\mathbb{L}^{\pm \Ind_\pi/2}(\Phi^0_{0,\pi})^\ast[d^0_\pi](p^\pm_V)_!(\eta^\pm_V)^\ast\\
    (p^\pm_{U_0})_!(\eta^\pm_{U_0})^\ast\phi_f\Phi^\ast[d]\arrow[r] & \phi_{f^0}(p^\pm_U)_!(\eta^\pm_U)^\ast\Phi^\ast[d]\arrow[ur,"\simeq"]
\end{tikzcd}\]
hence this isomorphism is compatible with pullbacks by smooth maps.\medskip

Consider two smooth spaces $U,V$ with a $\mathbb{C}^\ast$ action. The product space $U\times V$ is itself provided with a $\mathbb{C}^\ast$ action, and functors $X^0,X^\pm$ are obviously compatible with products, hence the hyperbolic localization diagram of $U\times V$ is the product of the hyperbolic localization diagrams of $U$ and $V$. In particular:
\begin{align}
    (p_{U\times V})_!(\eta_{U\times V})^\ast&=(p_U\times p_V)_!(\eta_U\times\eta_V)^\ast:D^b_c(U\times V)\to D^b_c(U^0\times V^0=(U\times V)^0)\nn\\
    (p_{U_0\times V_0})_!(\eta_{U_0\times V_0})^\ast&=(p_{U_0}\times p_{V_0})_!(\eta_{U_0}\times\eta_{V_0})^\ast:D^b_c(U_0\times V_0)\to D^b_c(U_0^0\times V_0^0=(U_0\times V_0)^0)
\end{align}
Denote now by $k:U_0\times V_0\to (U\times V)_0$ the natural closed embedding closed embedding. Using the two commutative diagram of \ref{diagcomthomseb}, one obtains that the following diagram is commutative:
\[\begin{tikzcd}[scale cd=0.8]\label{diagthomsebhyploc}
    (p^\pm_{U_0}\times p^\pm_{V_0})_!(\eta^\pm_{U_0}\times\eta^\pm_{V_0})^\ast(\phi_f\boxtimes\phi_g)\arrow[r]\arrow[d,"\simeq"] & (p^\pm_{U_0}\times p^\pm_{V_0})_!(\phi_{f^+}\boxtimes\phi_{g^+})(\eta^\pm_U\times\eta^\pm_V)^\ast\arrow[r]\arrow[d,"\simeq"]& (\phi_{f^0}\boxtimes\phi_{g^0})(p^\pm_U\times p^\pm_V)_!(\eta^\pm_U\times\eta^\pm_V)^\ast\arrow[d,"\simeq"]\\
    (k^0)^\ast(p^\pm_{U_0}\times p^\pm_{V_0})_!(\eta^\pm_{U_0}\times\eta^\pm_{V_0})^\ast\phi_{f\boxplus g}\arrow[r] & (k^0)^\ast(p^\pm_{U_0}\times p^\pm_{V_0})_!\phi_{f^+\boxplus g^+}(\eta^\pm_U\times\eta^\pm_V)^\ast\arrow[r] & (k_0)^\ast\phi_{f^0\boxplus g^0}(p^\pm_U\times p^\pm_V)_!(\eta^\pm_U\times\eta^\pm_V)^\ast
\end{tikzcd}\]
hence the isomorphism is compatible with Thom-Sebastiani isomorphism.\medskip

Again the same compatibility results can be lifted at the level of complexes of mixed Hodge modules.$\Box$\medskip

\subsection{Hyperbolic localization of the perverse sheaf of vanishing cycles}

Consider a critical chart $(R,U,f,i)$, \ie a smooth algebraic space $U$ with a regular function $f:U\to\mathbb{C}$ with $i:R\to U$ denoting the embedding of the critical locus. The image of $\phi_f$ is supported on $R$, hence denoting $i_c:R_c:=R\cap U_c\to U_c$ the embedding, $(i_c)^\ast\phi_f\simeq(i_c)^!\phi_f$ is a functor from $\Perv(U)$ to $\Perv(R_c)$. One defines in particular the perverse sheaf on $R$:
\begin{align}
    \mathcal{PV}_{U,f}:=\bigoplus_{c\in f(R)}i_c^\ast\phi_{f-c}\mathcal{IC}_U
\end{align}
with a monodromy operator $\tau_{U,f}:=\mathcal{PV}_{U,f}\to\mathcal{PV}_{U,f}$ defined by:
\begin{align}
    \tau_{U,f}=\bigoplus_{c\in f(R)}\tilde{T}_{f-c}^{2\pi}|_{R_c}
\end{align}
and a polarization operator $\sigma_{U,f}:\mathcal{PV}_{U,f}\overset{\simeq}{\to} \mathbb{D}\mathcal{PV}_{U,f}$ defined by the composition of isomorphisms:
\[\begin{tikzcd}
    \mathcal{PV}_{U,f}=\bigoplus_{c\in f(R)}i_c^\ast\phi_{f-c}\mathcal{IC}_U\arrow[r,"\simeq"] & \bigoplus_{c\in f(R)}i_c^\ast\phi_{f-c}\mathbb{D}\mathcal{IC}_V\arrow[d,"\simeq"]\\
    \mathbb{D}\mathcal{PV}_{U,f}=\mathbb{D}\bigoplus_{c\in f(R)}i_c^\ast\phi_{f-c}\mathcal{IC}_U & \bigoplus_{c\in f(R)}i_c^\ast\mathbb{D}\phi_{f-c}\mathcal{IC}_U\arrow[l,"\simeq"]
\end{tikzcd}\]
where we have used first the self-duality of the intersection complex, secondly Massey's isomorphism, and thirdly the fact that $(i_0)^\ast\phi_f=(i_0)^!\phi_f$.

This perverse sheaf is given the structure of strongly polarized a monodromic mixed Hodge module in \cite[sec 2.10]{Joycesymstab}, using the polarized mixed Hodge module on $\mathcal{IC}_U$, the strongly polarized monodromic mixed Hodge module obtained by applying $\phi^H_f$, and then using $i^\ast$ at the level of mixed Hodge modules. The functor $rat$ sends then the monodromy automorphism to $\tau_{U,f}$, and the strong polarization to $\sigma_{U,f}$.\medskip\medskip

\begin{proposition}\label{propchi}
    For a $\mathbb{C}^\ast$-equivariant critical chart $(U,R,f,i)$, there are natural isomorphisms in the abelian category $\bigoplus_{\pi\in\Pi}\MMHM(R^0_\pi)[-\Ind_\pi]$:
    \begin{align}
        \beta^\pm_{U,f}:(p_R^\pm)_!(\eta_R^\pm)^\ast\mathcal{PV}_{U,f}\to\bigoplus_{\pi\in\Pi}\mathbb{L}^{\pm \Ind_\pi/2}\mathcal{PV}_{U^0_\pi,f^0_\pi}
    \end{align}
    compatible with the strong polarization, \ie such that the following square is commutative:
    \[\begin{tikzcd}
    (p_R^+)_!(\eta_R^+)^\ast\mathcal{PV}_{U,f}\arrow[r,"\beta^+_{U,f}"]\arrow[d,"B_R^{-1}\circ\sigma_{U,f}"] & \bigoplus_{\pi\in\Pi}\mathbb{L}^{\Ind_\pi/2}\mathcal{PV}_{U^0_\pi,f^0_\pi}\arrow[d,"\bigoplus_{\pi\in\Pi}\sigma_{U^0_\pi,f^0_\pi}"]\\
    \mathbb{D}(p_R^-)_!(\eta_R^-)^\ast\mathcal{PV}_{U,f}\arrow[r,"\beta_{U,f}"] & \bigoplus_{\pi\in\Pi}\mathbb{L}^{-\Ind_\pi/2}\mathbb{D}\mathcal{PV}_{U^0_\pi,f^0_\pi}
    \end{tikzcd}\]
\end{proposition}

Proof: Consider the following commutative diagram:
\[\begin{tikzcd}
    U_c & U^\pm_c\arrow[l,"\eta^\pm_{U_c}"]\arrow[r,"p^\pm_{U_c}"] & U^0_c\\
    & (p^\pm_{U_c})^{-1}(R^0_c)\arrow[u,"\bar{i}"]\arrow[r,"\bar{p}"]& R^0_c\arrow[u,"i^0_c"]\\
    R_c\arrow[uu,"i_c"] & R_c^\pm\arrow[ur,swap,"p^\pm_R"] \arrow[l,"\eta^\pm_R"]\arrow[u,"\hat{i}"] &
\end{tikzcd}\]
where $\hat{i}$ and $\bar{i}$ are the obvious inclusions, and $\bar{p}:=p^\pm_{U_c}|_{(p^\pm_{U_c})^{-1}(R^0_c)}$. The upper right square is Cartesian. One can define the following isomorphisms:
\begin{align}\label{comihyp}
    (p^\pm_{R_c})_!(\eta^\pm_{R_c})^\ast (i_c)^\ast&\simeq\bar{p}_!\hat{i}_!\hat{i}^\ast\bar{i}^\ast(\eta^\pm_{U_c})^\ast\nn\\
    &\simeq\bar{p}_!\bar{i}^\ast(\eta_{U_c})^\ast\nn\\
    &\simeq(i^0_c)^\ast(p^\pm_{U_c})_!(\eta^\pm_{U_c})^\ast
\end{align}
where we have used the fact that $\hat{i}_!\hat{i}^\ast=\Id$ because $\hat{i}$ is a closed embedding in the second line, and the base change theorem in the upper right square in the last line. Moreover, because these isomorphism, and also the isomorphisms $S_R,S_{U_c}$ are built using the natural morphisms of the six-functor formalism in the above diagram, the following diagram is commutative, for $A\in D^b_c(U_c)_{R_c}^{\mathbb{C}^\ast-eq}$:
\[\begin{tikzcd}\label{diagcomidual}
    (p^+_{R_c})_!(\eta_{R_c}^+)^\ast i_c^\ast\mathbb{D}(A)\arrow[r,"\simeq"]\arrow[d,"\simeq"] & (i_c^0)^\ast(p^+_{U_c})_!(\eta^+_{U_c})^\ast\mathbb{D}(A)\arrow[dr,"\simeq"] &\\
    (p^+_{R_c})_!(\eta_{R_c}^+)^\ast\mathbb{D}i_c^\ast(A)\arrow[dr,"\simeq"] && (i_c^0)^\ast\mathbb{D}(p^-_{U_c})_!(\eta^-_{U_c})^\ast(A)\arrow[d,"\simeq"]\\
    & \mathbb{D}(p^-_{R_c})_!(\eta_{R_c}^-)^\ast i_c^\ast(A)\arrow[r,"\simeq"] & \mathbb{D}(i_c^0)^\ast(p^-_{U_c})_!(\eta^-_{U_c})^\ast(A)
\end{tikzcd}\]

Then one obtains isomorphisms:
\begin{align}
    \beta^\pm_{U,f}:(p^\pm_R)_!(\eta^\pm_R)^\ast\mathcal{PV}_{U,f}\simeq\bigoplus_{\pi\in\Pi}\mathbb{L}^{^\pm \Ind_\pi/2}{PV}_{U^0_\pi,f^0_\pi}
\end{align}
by the following composition of isomorphisms:
\[\begin{tikzcd}[scale cd=0.9]
    (p_R^\pm)_!(\eta_R^\pm)^\ast\mathcal{PV}_{U,f}:=(p_R^\pm)_!(\eta_R^\pm)^\ast\bigoplus_{c\in f(R)}i_c^\ast\phi_{f-c}\mathcal{IC}_U\arrow[r,"\simeq"] & \bigoplus_{c\in f(R)}(i_c^0)^\ast (p^\pm_{U_c})_!(\eta^\pm_{U_c})^\ast\phi_{f-c}\mathcal{IC}_U\arrow[d,"\simeq"]\\
    \bigoplus_{\pi\in\Pi}\mathbb{L}^{\pm \Ind_\pi/2}{PV}_{U^0_\pi,f^0_\pi}:=\bigoplus_{c\in f(R)}i_c^\ast\phi_{f-c}\bigoplus_{\pi\in\Pi}\mathbb{L}^{\pm \Ind_\pi/2}\mathcal{IC}_{U^0_\pi} & \bigoplus_{c\in f(R)}i_c^\ast\phi_{f-c}(p_U^\pm)_!(\eta_U^\pm)^\ast\mathcal{IC}_U\arrow[l,"\simeq"]
\end{tikzcd}\]
where the first isomorphism is \eqref{comihyp}, the second from proposition \ref{comhyplocvan}  where one has remaked that $\mathcal{IC}_U$ is $\mathbb{C}^\ast$-equivariant, and the third from corollary \ref{classbibi}.\medskip

Using the commutativity of the diagram \ref{diagfunctmon}, one obtains that the following square is commutative:
\[\begin{tikzcd}
    (p_{U_0}^\pm)_!(\eta_{U_0}^\pm)^\ast\phi_f\mathcal{IC}_U\arrow[r,"\simeq"]\arrow[d,"(p_{U_0}^\pm)_!(\eta_{U_0}^\pm)^\ast(\tilde{T}_f^{2\pi})"] & \bigoplus_{\pi\in\Pi}\mathbb{L}^{\pm \Ind_\pi/2}\phi_{f^0_\pi}\mathcal{IC}_{U^0_\pi}\arrow[d,"\bigoplus_{\pi\in\Pi}\tilde{T}^{2\pi}_{f^0_\pi}"]\\
    (p_{U_0}^\pm)_!(\eta_{U_0}^\pm)^\ast\phi_f\mathcal{IC}_U\arrow[r,"\simeq"] & \bigoplus_{\pi\in\Pi}\mathbb{L}^{\pm \Ind_\pi/2}\phi_{f^0_\pi}\mathcal{IC}_{U^0_\pi}
\end{tikzcd}\]
the compatibility with the monodromy results then by applying $(i_0)^\ast$ to this diagram, applying the functorial isomorphism \eqref{comihyp}, and finally writing the definition of $\tau_{U,f}$ by noticing that $(-1)^{\dim(U)}=(-1)^{\pm \Ind_\pi+\dim(U^0_\pi)}$.\medskip

The proof of the compatibility with $\sigma_{U,f}$ is slightly more technical. In fact $\beta^\pm_{U,f}$ and $\sigma_{U,f}$ are both defined by composing three isomorphisms expressing the commutation of the hyperbolic localization and Verdier duality with $i_c^\ast$, with $\phi_f$ and with $\mathcal{IC}_U$, and then one should divide the square expressing the compatibility into nine sub-squares:
\[\begin{tikzcd}
11 \arrow[r]\arrow[d] & 12 \arrow[r]\arrow[d] & 13 \arrow[r]\arrow[d] & 14 \arrow[d]\\
21 \arrow[r]\arrow[d] & 22 \arrow[r]\arrow[d] & 23 \arrow[r]\arrow[d] & 24 \arrow[d]\\
31 \arrow[r]\arrow[d] & 32 \arrow[r]\arrow[d] & 33 \arrow[r]\arrow[d] & 34 \arrow[d]\\
41 \arrow[r] & 42 \arrow[r] & 43 \arrow[r] & 44\\
\end{tikzcd}\]

Where the objects are:
\begin{align}
    &11:(p_R^+)_!(\eta_R^+)^\ast(i_0)^\ast\phi_f\mathcal{IC}_U\nn\\
    &12:(i^0_0)^\ast(p^+_{U_0})_!(\eta^+_{U_0})^\ast\phi_f\mathcal{IC}_U\nn\\
    &13:(i^0_0)^\ast\phi_{f^0}(p^+_U)_!(\eta^+_U)^\ast\mathcal{IC}_U\nn\\
    &14:(i^0_0)^\ast\phi_{f^0}\bigoplus_{\pi\in\Pi}\mathbb{L}^{\Ind_\pi/2}\mathcal{PV}_{U^0_\pi,f^0_\pi}\nn\\
    &21:(p_R^+)_!(\eta_R^+)^\ast(i_0)^\ast\phi_f\mathbb{D}\mathcal{IC}_U\nn\\
    &22:(i^0_0)^\ast(p^+_{U_0})_!(\eta^+_{U_0})^\ast\phi_f\mathbb{D}\mathcal{IC}_U\nn\\
    &23:(i^0_0)^\ast\phi_{f^0}(p^+_U)_!(\eta^+_U)^\ast\mathbb{D}\mathcal{IC}_U=(i^0_0)^\ast\phi_{f^0}\mathbb{D}(p^-_U)_!(\eta^-_U)^\ast\mathcal{IC}_U\nn\\
    &24:(i^0_0)^\ast\phi_{f^0}\bigoplus_{\pi\in\Pi}\mathbb{D}\mathbb{L}^{-\Ind_\pi/2}\mathcal{PV}_{U^0_\pi,f^0_\pi}\nn\\
    &31:(p_R^+)_!(\eta_R^+)^\ast(i_0)^\ast\mathbb{D}\phi_f\mathcal{IC}_U\nn\\
    &32:(i^0_0)^\ast(p^+_{U_0})_!(\eta^+_{U_0})^\ast\mathbb{D}\phi_f\mathcal{IC}_U=(i^0_0)^\ast\mathbb{D}(p^-_{U_0})_!(\eta^-_{U_0})^\ast\phi_f\mathcal{IC}_U\nn\\
    &33:(i^0_0)^\ast\phi_{f^0}(p^+_U)_!(\eta^+_U)^\ast\mathbb{D}\mathcal{IC}_U=(i^0_0)^\ast\mathbb{D}\phi_{f^0}(p^-_U)_!(\eta^-_U)^\ast\mathcal{IC}_U\nn\\
    &34:(i^0_0)^\ast\mathbb{D}\phi_{f^0}\bigoplus_{\pi\in\Pi}\mathbb{L}^{-\Ind_\pi/2}\mathcal{PV}_{U^0_\pi,f^0_\pi}\nn\\
    &41:(p_R^+)_!(\eta_R^+)^\ast\mathbb{D}(i_0)^\ast\phi_f\mathcal{IC}_U=\mathbb{D}(p_R^-)_!(\eta_R^-)^\ast(i_0)^\ast\phi_f\mathcal{IC}_U\nn\\
    &42:\mathbb{D}(i^0_0)^\ast(p^+_{U_0})_!(\eta^+_{U_0})^\ast\phi_f\mathcal{IC}_U=(i^0_0)^\ast\mathbb{D}(p^-_{U_0})_!(\eta^-_{U_0})^\ast\phi_f\mathcal{IC}_U\nn\\
    &43:\mathbb{D}(i^0_0)^\ast\phi_{f^0}(p^+_U)_!(\eta^+_U)^\ast\mathcal{IC}_U=(i^0_0)^\ast\mathbb{D}\phi_{f^0}(p^-_U)_!(\eta^-_U)^\ast\mathcal{IC}_U\nn\\
    &44:\mathbb{D}(i^0_0)^\ast\phi_{f^0}\bigoplus_{\pi\in\Pi}\mathbb{L}^{-\Ind_\pi/2}\mathcal{PV}_{U^0_\pi,f^0_\pi}
\end{align}
The down left square commutes because it is the diagram \ref{diagcomidual} with $A=\phi_f\mathcal{IC}_U$. The central square commutes because \ref{diagcomvandual} commutes, and the upper right square commutes because it is $(i^0_0)^\ast\phi_{f^0}$ applied to the commutative diagram in corollary \ref{classbibi}. The squares which are not off the down-left/up right diagonal are commutative because the isomorphisms constructed are functorial. Then the whole diagram is commutative, and, summing over $c\in f(R)$, one obtains the commutativity of the diagram expressing the compatibility with the polarization.\medskip

Again the same compatibility results can be directly lifted to the derived category of monodromic mixed Hodge modules.$\Box$\medskip

\subsection{Compatibility with smooth maps}

Consider two critical charts $(R,U,f,i),(S,V,g,j)$ and $\Phi:U\to V$ a smooth map of relative dimension $d$ with $f=g\circ\Phi$. Joyce defines then in \cite{darbstack} the isomorphism
\begin{align}
    \Xi_\Phi:(\Phi|_R)^\ast[d]\mathcal{PV}_{V,g}\simeq\mathcal{PV}_{U,f}
\end{align}
by the following composition of isomorphisms:
\[\begin{tikzcd}
    (\Phi|_R)^\ast[d]\mathcal{PV}_{V,g}=(\Phi|_R)^\ast[d]\bigoplus_{c\in g(S)}j_c^\ast\phi_{g-c}\mathcal{IC}_V\arrow[r,"\simeq"] & \bigoplus_{c\in f(R)}i_c^\ast (\Phi_c)^\ast[d]\phi_{g-c}\mathcal{IC}_V\arrow[d,"\simeq"]\\
    \mathcal{PV}_{U,f}=\bigoplus_{c\in f(R)}i_c^\ast\phi_{f-c}\mathcal{IC}_U & \bigoplus_{c\in f(R)}i_c^\ast\phi_{f-c}\Phi^\ast[d]\mathcal{IC}_V\arrow[l,"\simeq"]
\end{tikzcd}\]
identifying $(\Phi|_R)^\ast[d](\sigma_{V,g})$ with $\sigma_{U,f}$ and $(\Phi|_R)^\ast[d](\tau_{V,g})$ with $\tau_{U,f}$ according to proposition 4.3 of \cite{darbstack}.\medskip

\begin{proposition}\label{propcomsmooth}
    The isomorphisms $\beta^\pm_{U,f}$ commute with pullback of smooth maps of $\mathbb{C}^\ast$-equivariant critical charts, \ie for $\Phi:U\to V$ a $\mathbb{C}^\ast$ equivariant and smooth map of relative dimension $d$ and critical charts $(R,U,f,i)$ and $(S,V,g,j)$ with $f=g\circ\Phi$ one has the commutative diagram:
    \[\begin{tikzcd}[row sep=5ex,column sep=5ex]\label{diagcomsmooth}
        \bigoplus_{\pi\in\Pi}\mathbb{L}^{\pm \Ind_\pi/2}(\Phi^0_\pi|_{R^0})^\ast[d^0_\pi](p^\pm_S)_!(\eta^\pm_S)^\ast\mathcal{PV}_{V,g}\arrow[r,"\beta^\pm_{V,g}"]  & \bigoplus_{\pi,\pi'\in\Pi}\mathbb{L}^{(\pm \Ind_{\pi'})/2}(\Phi^0_\pi|_{R^0})^\ast[d^0_\pi]\mathcal{PV}_{V^0_{\pi'-\pi},g^0_{\pi'-\pi}}\arrow[dd,"\Xi_{\Phi^0}"]\\
        (p^\pm_R)_!(\eta^\pm_R)^\ast(\Phi|_R)^\ast[d]\mathcal{PV}_{V,g}\arrow[d,"\Xi_\Phi"]\arrow[u,"\simeq"]&\\
        (p^\pm_R)_!(\eta^\pm_R)^\ast\mathcal{PV}_{U,f}\arrow[r,"\beta^\pm_{U,f}"] & \bigoplus_{\pi\in\Pi}\mathbb{L}^{\pm \Ind_{\pi'}/2}\mathcal{PV}_{U^0_ {\pi'},f^0_{\pi'}}
    \end{tikzcd}\]
\end{proposition}

In particular, one obtains that $\beta^\pm_{U,f}$ commute with the restriction on an étale neighborhood.

Proof: Consider now the following diagram, where we denote $G=\phi_g\mathcal{IC}_V$:
\[\begin{tikzcd}[column sep=-8ex]\label{diag1}
    (p_R)_!(\eta_R)^\ast(\Phi|_R)^\ast(j_0)^\ast G=(\Phi|_{R^0})^\ast(p_S)_!(\eta_S)^\ast(j_0)^\ast G\arrow[r]\arrow[d] & (\Phi|_{R^0})^\ast(j^0_0)^\ast(p_{V_0})_!(\eta_{V_0})^\ast G\arrow[d]\\
    (p_R)_!(\eta_R)^\ast(i_0)^\ast(\Phi_0)^\ast G\arrow[r] &(i^0_0)^\ast(p_{U_0})_!(\eta_{U_0})^\ast (\Phi_0)^\ast G=(i^0_0)^\ast(\Phi^0_0)^\ast(p_{V_0})_!(\eta_{V_0})^\ast
\end{tikzcd}\]
where the horizontal isomorphisms are those of \eqref{comihyp} applied to $V$ and $U$: it is a diagram of isomorphisms between objects in an abelian category, and the isomorphisms are built from the natural isomorphisms of the six-functor formalism in the diagram:
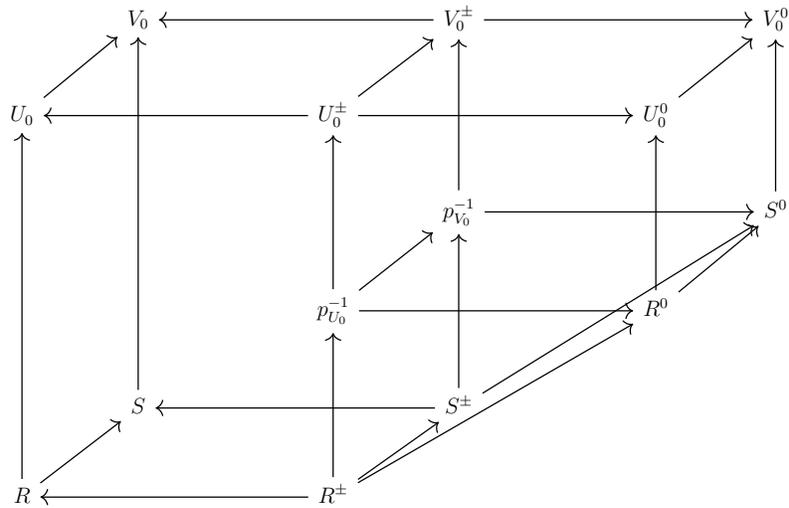
\[\begin{tikzcd}[scale cd=0.7]
&V_0 &&& V_0^\pm\arrow[lll]\arrow[rrr] &&& V_0^0\\
U_0\arrow[ur] &&& U_0^\pm\arrow[ur]\arrow[lll]\arrow[rrr] &&& U_0^0\arrow[ur]&\\
&  &&& p^{-1}_{V_0}\arrow[uu]\arrow[rrr] &&& S^0\arrow[uu] \\
&&& p^{-1}_{U_0}\arrow[uu]\arrow[rrr]\arrow[ur] &&& R^0\arrow[uu]\arrow[ur] &\\
&S \arrow[uuuu] &&& S^\pm\arrow[uu]\arrow[uurrr]\arrow[lll] &&& \\
R \arrow[uuuu]\arrow[ur] &&& R^\pm\arrow[uu]\arrow[uurrr]\arrow[ur]\arrow[lll] &&&&\\
\end{tikzcd}\]

hence this diagram commutes. Finally, consider the diagram:
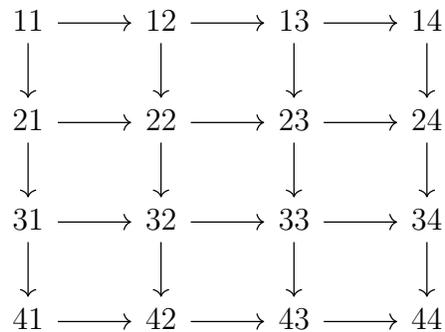
\[\begin{tikzcd}
11 \arrow[r]\arrow[d] & 12 \arrow[r]\arrow[d] & 13 \arrow[r]\arrow[d] & 14 \arrow[d]\\
21 \arrow[r]\arrow[d] & 22 \arrow[r]\arrow[d] & 23 \arrow[r]\arrow[d] & 24 \arrow[d]\\
31 \arrow[r]\arrow[d] & 32 \arrow[r]\arrow[d] & 33 \arrow[r]\arrow[d] & 34 \arrow[d]\\
41 \arrow[r] & 42 \arrow[r] & 43 \arrow[r] & 44
\end{tikzcd}\]

Where the objects are:
\begin{align}
    &11:(p_R)_!(\eta_R)^\ast(\Phi|_R)^\ast[d](j_0)^\ast\phi_g\mathcal{IC}_V=(\Phi|_{R^0})^\ast[d](p_S)_!(\eta_S)^\ast(j_0)^\ast\phi_g\mathcal{IC}_V\nn\\
    &12:(\Phi|_{R^0})^\ast[d](j^0_0)^\ast(p_{V_0})_!(\eta_{V_0})^\ast\phi_g\mathcal{IC}_V\nn\\
    &13:(\Phi|_{R^0})^\ast[d](j^0_0)^\ast\phi_{g^0}(p_V)_!(\eta_V)^\ast\mathcal{IC}_V\nn\\
    &14:(\Phi|_{R^0})^\ast[d](j^0_0)^\ast\phi_{g^0}\bigoplus_{n\in\mathbb{Z}}\mathbb{L}^{n/2}\mathcal{PV}_{V^0_n,g^0_n}\nn\\
    &21:(p_R)_!(\eta_R)^\ast(i_0)^\ast(\Phi_0)^\ast[d]\phi_g\mathcal{IC}_V\nn\\
    &22:(i^0_0)^\ast(p_{U_0})_!(\eta_{U_0})^\ast(\Phi_0)^\ast[d]\phi_g\mathcal{IC}_V=(i^0_0)^\ast(\Phi^0_0)^\ast[d](p_{V_0})_!(\eta_{V_0})^\ast\phi_g\mathcal{IC}_V\nn\\
    &23:(i^0_0)^\ast(\Phi^0_0)^\ast[d]\phi_{g^0}(p_V)_!(\eta_V)^\ast\mathcal{IC}_V\nn\\
    &24:(i^0_0)^\ast(\Phi^0_0)^\ast[d]\phi_{g^0}\bigoplus_{n\in\mathbb{Z}}\mathbb{L}^{n/2}\mathcal{PV}_{V^0_n,g^0_n}\nn\\
    &31:(p_R)_!(\eta_R)^\ast(i_0)^\ast\phi_f\Phi^\ast[d]\mathcal{IC}_V\nn\\
    &32:(i^0_0)^\ast(p_{U_0})_!(\eta_{U_0})^\ast\phi_f\Phi^\ast[d]\mathcal{IC}_V\nn\\
    &33:(i^0_0)^\ast\phi_{f^0}(p_U)_!(\eta_U)^\ast\Phi^\ast[d]\mathcal{IC}_V=(i^0_0)^\ast\phi_{f^0}(\Phi^0)^\ast[d](p_V)_!(\eta_V)^\ast\mathcal{IC}_V\nn\\
    &34:(i^0_0)^\ast\phi_{f^0}(\Phi^0)^\ast[d]\bigoplus_{n\in\mathbb{Z}}\mathbb{L}^{n/2}\mathcal{PV}_{V^0_n,g^0_n}\nn\\
    &41:(p_R)_!(\eta_R)^\ast(i_0)^\ast\phi_f\mathcal{IC}_U\nn\\
    &42:(i^0_0)^\ast(p_{U_0})_!(\eta_{U_0})^\ast\phi_f\mathcal{IC}_U\nn\\
    &43:(i^0_0)^\ast\phi_{f^0}(p_U)_!(\eta_U)^\ast\mathcal{IC}_U\nn\\
    &44:(i^0_0)^\ast\phi_{f^0}\bigoplus_{n\in\mathbb{Z}}\mathbb{L}^{n/2}\mathcal{PV}_{U^0_n,f^0_n}
\end{align}
the three diagonal squares of the above diagram are respectively the commutative diagrams \ref{diag1}, \ref{diagcomsmoothvanhyp} and the commutative diagram expressing the fact that the isomorphism giving the commutation of the hyperbolic localization and pullback by smooth maps in proposition \ref{propsmoothhyploc} is compatible with composition of smooth maps, and the off diagonal squares commute because the functoriality of the isomorphism. Then the whole diagram is commutative. According to the definition of $\Xi_\Phi,\Xi_{\Phi^0}$ and $\beta_{U,f},\beta_{V,g}$, the diagram \ref{diagcomsmooth} is the direct sum of the above diagram with $f$ replaced by $f-c$ for $c\in f(R)$, hence it is also commutative. $\Box$

\subsection{Compatibility with Thom-Sebastiani isomorphism}

Consider two critical charts $(U,R,f,i)$ and $(S,V,g,j)$. Notice that because $U$ and $V$ are smooth $\Crit(f\boxplus g)=\Crit(f)\times \Crit(g)R\times S$. Consider the sequence of closed embeddings:
\[\begin{tikzcd}
R\times S\arrow[r,"i_c\times j_d"] &U_c\times V_d\arrow[r,"k_{c,d}"]& (U\times V)_{c+d}
\end{tikzcd}\]
In \cite{Joycesymstab}, Joyce defined the Thom-Sebastiani isomorphism:
\begin{align}
    \mathcal{TS}_{U,f,V,g}:\mathcal{PV}_{U\times V,f\boxplus g}\simeq\mathcal{PV}_{U,f}\boxtimes\mathcal{PV}_{S,g}
\end{align}
By the following sequence of isomorphisms:
\begin{align}
    \mathcal{PV}_{U\times V,f\boxplus g}=&\bigoplus_{e\in f\boxplus g(R\times S)}((i\times j)_e)^\ast\phi_{f\boxplus g-e}\mathcal{IC}_{U\times V}\nn\\
    \simeq &\bigoplus_{c\in f(R),d\in g(S)}(i_c\times j_d)^\ast k_{c,d}^\ast\phi_{(f-c)\boxplus(g-d)}\mathcal{IC}_{U\times V}\nn\\
    \simeq&\bigoplus_{c\in f(R),d\in g(S)}(i_c\times j_d)^\ast(\phi_{f-c}\boxtimes\phi_{g-d})\mathcal{IC}_{U\times V}\nn\\
    \simeq&(\bigoplus_{c\in f(R)}(i_c)^\ast\phi_{f-c}\mathcal{IC}_{U})\boxtimes(\bigoplus_{d\in g(S)}(j_d)^\ast\phi_{g-d}\mathcal{IC}_{V})\nn\\
    =&\mathcal{PV}_{U,f}\boxtimes\mathcal{PV}_{S,g}
\end{align}
where the third line is the Thom-Sebastiani isomorphism of \cite{MasThomSeb}.\medskip
    
\begin{proposition}\label{propcomthomseb}
    The isomorphism $\beta_{U,f}$ commutes with Thom-Sebastiani isomorphism, namely for $\mathbb{C}^\ast$-equivariant critical chart $(U,R,f,i)$ and $(S,V,g,j)$, one has the commutative diagram:
    \[\begin{tikzcd}[row sep=10ex,column sep=12ex]
    \scriptsize\parbox{6cm}{$(p_R^\pm\times p_S^\pm)|_!(\eta_R^\pm\times\eta_S^\pm)^\ast(\mathcal{PV}_{U\times V,f\boxplus g})$}\arrow[r,"\beta^\pm_{U\times V}"]\arrow[d,"\scriptsize\parbox{4cm}{$(p_R^\pm\times p_S^\pm)|_!(\eta_R^\pm\times\eta_S^\pm)^\ast\mathcal{TS}_{U,f,V,g}$}"] & 
    \scriptsize\parbox{6cm}{$\bigoplus_{\pi,\pi'\in\Pi}\mathbb{L}^{\pm(\Ind_\pi+\Ind_{\pi'})/2}\mathcal{PV}_{U^0_{\pi}\times V^0_{\pi'},f^0_{\pi}\boxplus g^0_{\pi'}}$}\arrow[d,swap,"\scriptsize\parbox{4cm}{$\bigoplus_{\pi,\pi'\in\Pi}\mathbb{L}^{\pm(\Ind_\pi+\Ind_{\pi'})/2}\mathcal{TS}_{U^0_{\pi},f^0_{\pi},V^0_{\pi'},g^0_{\pi'}}$}"]\\
    \scriptsize\parbox{6cm}{$((p_R^\pm)_!(\eta_R^\pm)^\ast\mathcal{PV}_{U,f})\boxtimes((p_S^\pm)_!(\eta_S^\pm)^\ast\mathcal{PV}_{S,g})$}\arrow[r,"\beta_U\times\beta_V"] &
    \scriptsize\parbox{6cm}{$\bigoplus_{\pi\in\mathbb{Z}}\mathbb{L}^{\pm \Ind_\pi/2}\mathcal{PV}_{U^0_{\pi},f^0_{\pi}})\boxtimes\bigoplus_{\pi'\in\Pi}\mathbb{L}^{\pm \Ind_{\pi'}/2}\mathcal{PV}_{V^0_{\pi'},g^0_{\pi'}})$}
    \end{tikzcd}\]
\end{proposition}

Proof: Consider the following diagram:
\[\begin{tikzcd}
11 \arrow[r]\arrow[d] & 12 \arrow[r]\arrow[d] & 13 \arrow[r]\arrow[d] & 14 \arrow[d]\\
21 \arrow[r] & 22 \arrow[r] & 23 \arrow[r] & 24 
\end{tikzcd}\]
Where the objects are:
\begin{align}
    &11:(p_R^\pm\times P_S^\pm)_!(\eta_R^\pm\times\eta_S^\pm)^\ast(i_0\times j_0)^\ast k^\ast\phi_{f\boxplus g}\mathcal{IC}_{U\times V}\nn\\
    &12:(i^0_0\times j^0_0)^\ast(p^\pm_{U_0}\times p^\pm_{V_0})_!(\eta^\pm_{U_0}\eta^\pm_{V_0})^\ast k^\ast\phi_{f\boxplus g}\mathcal{IC}_{U\times V}=(i^0_0\times j^0_0)^\ast(k^0)^\ast(p^\pm_{(U\times V)_0})_!(\eta^\pm_{(U\times V)_0})^\ast \phi_{f\boxplus g}\mathcal{IC}_{U\times V}\nn\\
    &13:(i^0_0\times j^0_0)^\ast(k^0)^\ast\phi_{f^0\boxplus g^0}(p^\pm_U\times p^\pm_V)_!(\eta^\pm_U\times\eta^\pm_V)^\ast \mathcal{IC}_{U\times V}\nn\\
    &14:(i^0_0\times j^0_0)^\ast(k^0)^\ast\phi_{f^0\boxplus g^0}\bigoplus_{\pi,\pi'\in\Pi}\mathbb{L}^{(\Ind_\pi+\Ind_{\pi'})/2}\mathcal{IC}_{U^0_{\pi}\times V^0_{\pi'}}\nn\\
    &21:(p_R^\pm\times P_S^\pm)_!(\eta_R^\pm\times\eta_S^\pm)^\ast(i_0\times j_0)^\ast (\phi_f\boxtimes\phi_g)\mathcal{IC}_{U\times V}\nn\\
    &22:(i^0_0\times j^0_0)^\ast(p^\pm_{U_0}\times p^\pm_{V_0})_!(\eta^\pm_{U_0}\eta^\pm_{V_0})^\ast (\phi_f\boxtimes\phi_g)\mathcal{IC}_{U\times V}\nn\\
    &23:(i^0_0\times j^0_0)^\ast(\phi_{f^0}\boxtimes\phi_{g^0})(p^\pm_U\times p^\pm_V)_!(\eta^\pm_U\times\eta^\pm_V)^\ast \mathcal{IC}_{U\times V}\nn\\
    &24:(i^0_0\times j^0_0)^\ast(\phi_{f^0}\boxtimes\phi_{g^0})\bigoplus_{\pi,\pi'\in\Pi}\mathbb{L}^{(\Ind_\pi+\Ind_{\pi'})/2}\mathcal{IC}_{U^0_{\pi}\times V^0_{\pi'}}
\end{align}
the left and the right squares commutes from the functoriality of the corresponding isomorphisms, and the central square commutes because it is the diagram \ref{diagthomsebhyploc}, which commutes, hence the whole diagram commutes. Identifying $\mathcal{IC}_{U\times V}$ with $\mathcal{IC}_U\boxtimes\mathcal{IC}_V$ and $\mathcal{IC}_{U^0_\pi\times V^0_{\pi'}}$ with $\mathcal{IC}_{U^0_\pi}\boxtimes\mathcal{IC}_{V^0_{\pi'}}$, and summing over $c\in f(R),d\in g(S)$, one obtains the commutative square of the proposition. $\Box$\medskip

As explained in example 2.14 of \cite{Joycesymstab}, for $q$ a non-degenerate quadratic form on an $n$-dimensional vector space $E$, one has $\Crit(q)=\{0\}$, and then:
\begin{align}
    \mathcal{PV}_{E,q}= H^{n-1}(MF_q(0),\mathbb{Q})\otimes \mathbb{Q}_{\{0\}}\simeq \mathbb{Q}_{\{0\}}
\end{align}
where $MF_q(0)$ denotes the Milnor fiber of $q$ at $0$, which is $T^\ast S^{n-1}$, and the second isomorphism comes from the orientation of $S^{n-1}$ coming from an orientation of $E$. Consider now $E$ with a linear $\mathbb{C}^\ast$-action and a $\mathbb{C}^\ast$-invariant non-degenerate quadratic form $q$.  We can decompose $E=E^0\oplus E^+\oplus E^-$ according to the $\mathbb{C}^\ast$-weights, and the non-degenerate invariant quadratic form $q$ gives a natural isomorphism $E^-=(E^+)^\vee$, hence a natural isomorphism $K_E=K_{E^0}$, \ie a natural bijection between orientations of $E$ and orientations of $E^0$. Denoting $q^0:=q|_{E^0}$, one obtains directly the commutativity of the following square:
\[\begin{tikzcd}\label{diagnatident}
\mathcal{PV}_{E,q}\arrow[r,"\beta^\pm_{E,q}"]\arrow[d,"\simeq"] & \mathcal{PV}_{E^0,q^0}\arrow[d,"\simeq"]\nn\\
\mathbb{Q}_{\{0\}}\arrow[r,"\simeq"] & \mathbb{Q}_{\{0\}}
\end{tikzcd}\]
where the vertical arrows comes from consistent orientations of $V$ and $V^0$.

\section{Białynicki-Birula decomposition on a d-critical algebraic space}

\subsection{D-critical algebraic space}

We recall here the notions and results of \cite{Joyce2013ACM} and \cite{Joycesymstab} about d-critical structures and the Donaldson-Thomas sheaf. Joyce and his collaborators have developed the theory of d-critical schemes and build the perverse sheaf on it using the Zariski topology, mainly in order to define also motivic Donaldson-Thomas invariants, which can be glued in the Zariski topology, but not in the étale topology. But because coherent sheaves and perverse sheaves glue in the étale topology, and the perverse sheaf of vanishing cycles transforms naturally under étale maps, the same formalism can be developed in the étale topology. We will then work in the étale topology, and consider d-critical algebraic spaces.\medskip

In \cite[Theo 2.1]{Joyce2013ACM}, Joyce constructs sheaves $\mathcal{S}_X$ and $\mathcal{S}^0_X$ on any algebraic space $X$. Given an étale map $R\to X$ and a closed embedding $i:R\hookrightarrow U$ into a smooth scheme $U$, denote by $I_{R,U}$ the sheaf of ideal in $i^{-1}(\mathcal{O}_U)$ of functions on $U$ near $i(R)$ which vanishes on $i(R)$. There is then an exact sequence of sheaves on $R$ defining locally $\mathcal{S}_X$:
\begin{align}
    0\to\mathcal{S}_X|_R\overset{\iota_{R,U}}{\to}\frac{i^{-1}(\mathcal{O}_U)}{I^2_{R,U}}\overset{d}{\to}\frac{i^{-1}(T^\ast U)}{I_{R,U}.i^{-1}(T^\ast U)}
\end{align}

One has a decomposition $\mathcal{S}_X=\mathcal{S}^0_X\oplus\mathbb{C}$, and $\mathcal{S}^0_X\subset\mathcal{S}_X$ is the kernel of the composition:
\begin{align}
    \mathcal{S}_X\to\mathcal{O}_X\to\mathcal{O}_{X^{red}}
\end{align}
with $X^{red}$ the reduced subspace of $X$. This construction is functorial, \ie for $\Phi: X\to Y$ a morphism of scheme there is a sheaf morphism $\Phi^\star:\mathcal{S}^0_Y\to\mathcal{S}^0_X$.\medskip

A d-critical structure on $X$ is a section $s$ of $\mathcal{S}^0_X$ such that for each $x\in X$ there exist an étale neighborhood $R$ of $x$, and embedding $i:R\hookrightarrow U$ into a smooth algebraic space $U$, and a regular function $f:U\to\mathbb{C}$ such that $i(R)=\Crit(f)$ and $\iota_{R,U}(s|_R)=i^{-1}(f)+I_{R,U}^2$. Informally, the data of $s$ precise the functional $f$ of critical charts up to second order terms, and one has $f|_{R^{red}}=0$. We deal then with critical charts $(R,U,f,i)$: considering an étale map $U'\to U$, one can consider a subchart $(R':=R\times_U U',U',f':=f|_{U'},i':=i|_{R'})\to(R,U,f,i)$.\medskip

One can also consider embeddings of charts $\Phi:(R,U,f,i)\hookrightarrow(S,V,g,j)$ for $R\to S\to X$ étale maps, \ie a locally closed embedding $\Phi:U\to V$ such that $\Phi\circ i=j|_R:R\to V$ and $f=g\circ\Phi:U\to\mathbb{C}$. According to \cite[Theo 2.20]{Joyce2013ACM}, one can compare two critical charts using embeddings: namely, for $(R,U,f,i)$ and $(S,V,g,j)$ two critical charts and $x\in R\cap S$, there exist subcharts $(R',U',f',i')\subset(R,U,f,i)$ and $(S',V',g',j')\subset(S,V,g,j)$ such that $x\in R'\cap S'$, a critical chart $(T,W,h,k)$ and embeddings $\Phi:(R',U',f',i')\hookrightarrow(T,W,h,k)$ and $\Psi(S',V',g',j')\hookrightarrow(T,W,h,k)$.\medskip

For a $-1$-shifted symplectic structure, on a d-critical chart $(R,U,f,i)$, the tangent-obstruction complex $\mathbb{L}_X|_{R^{red}}$ is quasi-isomorphic with $0\to TU\to T^\ast U\to 0$, and then $det(\mathbb{L}_X)|_{R^{red}}=i^\ast(K_U^{\otimes 2})|_{R^{red}}$, \ie the sheaves $i^\ast(K_U^{\otimes 2})|_{R^{red}}$ glue in a sheaf on $X^{red}$ For a general d-critical scheme and an embedding $\Phi:(R,U,f,i)\to(S,V,g,j)$, from definition 2.26 of \cite{Joyce2013ACM}, one has a natural isomorphism:
\begin{align}
    J_\Phi:i^\ast(K_U^{\otimes^2})|_{R^{red}}\to j^\ast(K^{\otimes^2}_V)|_{R^{red}}
\end{align}
Using the fact that two maps can be locally embedded in a single map, it is shown in \cite[Theo 2.28]{Joyce2013ACM} that these sheaves glue together into a single sheaf $K_{X,s}$ on $X^{red}$ called the canonical sheaf, with natural local isomorphisms:
\begin{align}
    &\iota_{R,U,f,i}:K_{X,s}|_{R^{red}}\to i^\ast(K_U^{\otimes^2})|_{R^{red}}\nn\\
    &J_\Phi\circ\iota_{R,U,f,i}=\iota_{S,V,g,j}\label{iotaphi}
\end{align}

\subsection{Gluing the Donaldson-Thomas sheaves}

Given a d-critical scheme, one could try naively to define a perverse sheaf and strongly polarized monodromic mixed Hodge modules modeled locally on $\mathcal{PV}_{U,f}$ for each critical charts $(R,U,f,i)$, gluing these sheaves and strongly polarized monodromic mixed Hodge modules by constructing isomorphisms on intersections of critical charts, satisfying the cocycle conditions. This is complicated by the orientations issue seen above. The construction of $\mathcal{PV}_{U,f}$ is natural with respect to étale restriction of maps, and étale locally two intersecting critical charts are related by stabilization \ie embedding of the form $\Phi:(R,U,f,i)\hookrightarrow(S,U\times E,g\boxplus q,j\times 0)$, with $E$ a vector space and $q$ a non-degenerate quadratic form. As seen above, the descent data to define the orientation $K_{X,s}^{1/2}$ is equivalent to give a natural orientation on $E$ for any such stabilization. We can then consider the chain of isomorphisms in $\Perv(X)$ (or $MMHM(X)$):
\begin{align}
    \Theta(\Phi):\mathcal{PV}_{U\times E,f\boxplus q}\simeq \mathcal{PV}_{U,f}\boxtimes\mathcal{PV}_{E,q}\simeq \mathcal{PV}_{U,f}
\end{align}
Where the first isomorphism is the Thom-Sebastiani isomorphism $\mathcal{TS}_{U,f,E,q}$, and the second comes from the natural orientation of $E$. The technical work of \cite{Joycesymstab} is to check that $(\mathcal{PV}_{U,f},\Theta(\Phi))$ defines a descent data, namely that the $\Theta(\Phi)$ glue on intersections of critical charts to define comparison isomorphisms, and that the cocycle relations are verified.\medskip

We are here interested in building an isomorphism between perverse sheaves, so we have to work one categorical level below: namely, we must define these isomorphisms locally, and check that this isomorphism commutes with the gluing isomorphisms. We will then use only the above presentation $(\mathcal{PV}_{U,f},\Theta(\Phi))$ for $P_{X,s}$, and will not use the equivalent, but more technical, presentation of \cite{Joycesymstab}.

\subsection{Białynicki-Birula decomposition}

An action $\mu:T\times X\to X$ of a one dimensional torus $\mathbb{C}^\ast$ on a d-critical scheme $X$ is said to leave invariant the d-critical structure $s$ (resp. the orientation $K_{X,s}^{1/2}$) if $\mu(\gamma)^\star(s)=s$ (resp. $\mu(\gamma)^{\ast}K_{X,s}^{1/2}=K_{X,s}^{1/2}$) for $\gamma\in T$. In particular, the d-critical structure of the classical truncation of a $(-1)$-shifted symplectic scheme with a $\mathbb{C}^\ast$-action leaving the $(-1)$-shifted symplectic structure invariant is $\mathbb{C}^\ast$-invariant. If the action is étale locally linearizable, then translating \cite[Prop 2.43, 2.44]{Joyce2013ACM} from the Zariski to the étale topology, we can then work with $\mathbb{C}^\ast$-equivariant critical charts. Namely, we can cover $X$ in the étale topology by charts $(R,U,f,i)$ such that $U$ has a $\mathbb{C}^\ast$-action for which $i$ is equivariant and $f$ is invariant. Moreover, considering $(R,U,f,i)$ and $(S,V,g,j)$  two $\mathbb{C}^\ast$-equivariant critical charts, and $x\in R\cap S$, one has étale restrictions $(R',U',f',i')\to(R,U,f,i)$ and $(S',V',g',j')\to(S,V,g,j)$ with $x\in R'\cap S'$, and $\mathbb{C}^\ast$-equivariant critical chart $(T,W,h,k)$ and $\mathbb{C}^\ast$-equivariant embeddings $\Phi:(R',U',f',i')\hookrightarrow(T,W,h,k)$ and $\Psi(S',V',g',j')\hookrightarrow(T,W,h,k)$.\medskip

As explained in \cite[Cor 2.45]{Joyce2013ACM}, considering the closed embedding $\xi_X:X^0\hookrightarrow X$, $X^0$ admits a natural d-critical structure $(X^0,s^0:=\xi_X^\star(s))$. Indeed, for $(R,U,f,i)$ a $\mathbb{C}^\ast$-equivariant d-critical chart, denoting by $R^0$, $U^0$ the $\mathbb{C}^\ast$-fixed locus of $R$,$U$, and $f^0=f|_{U^0}$, $i^0=i|_{R^0}$, $(R^0,U^0,f^0,i^0)$ is a d-critical chart of $X^0$ such that $s^0_{R^0,U^0}=f^0+I^2_{R^0,U^0}$, and then $X^0$ is covered by d-critical charts.\medskip

Consider now the orientation issue. Suppose first that $(X,s)$ is the classical truncation of a $-1$-shifted symplectic space with a $\mathbb{C}^\ast$-action. Then the tangent-obstruction complex $\mathbb{L}_X|_{X^0}$ splits as a direct sum:
\begin{align}
\mathbb{L}_X|_{X^0}=\mathbb{L}_{X^0}\oplus\mathbb{L}^+_{X^0}\oplus\mathbb{L}^-_{X^0}
\end{align}
Here $\mathbb{L}^+_{X^0}$ (resp $\mathbb{L}^-_{X^0}$) denotes the part of contracting (resp repelling) weight under the $\mathbb{C}^\ast$-action. The $-1$-shifted symplectic structure provides then a canonical isomorphism between $\mathbb{L}^+_{X^0}$ and $(\mathbb{L}^-_{X^0})^\vee[1]$, hence, taking:
\begin{align}
&det(\mathbb{L}_X)|_{X^0}=det(\mathbb{L}_{X^0})\otimes det(\mathbb{L}^+_{X^0})^2\nn\\
\implies & K_{X^0,s^0}=K_{X,s}\otimes det(\mathbb{L}^+_{X^0})^{-2}
\end{align}
Hence, given an orientation $K_{X,s}^{1/2}$ on $X$, we can define canonically an orientation on $X^0$:
\begin{align}
    K_{X^0,s^0}^{1/2}:=K_{X,s}^{1/2}\otimes det(\mathbb{L}^+_{X^0})^{-1}
\end{align}
Such construction is still possible if $(X,s)$ is a d critical space with $\mathbb{C}^\star$ action, not necessarily coming from a $-1$-shifted symplectic space:

\begin{lemma}\label{lemorient}
    If $(X,s)$ is a d critical space with $\mathbb{C}^\star$ action, there is a natural bijection between orientations of $(X,s)$ and orientations of $(X^0,s^0)$.
\end{lemma}

Proof: We can define an orientation $K_{X,s}^{1/2}$ (resp $K_{X^0,s^0}^{1/2}$) by descent: we define $(K_{X,s}^{1/2})_{R}:=K_{U}$ (resp $(K_{X^0,s^0}^{1/2})_{R^0}:=K_{U^0}$) for each critical chart $(R,U,f,i)$ (resp $(R^0,U^0,f^0,i^0)$, and gluing isomorphisms on double intersections which satisfies cocycle on triple intersections. Using \cite[Prop 2.43]{Joyce2013ACM}, we can use $\mathbb{C}^\ast$-equivariant critical charts, and, using \cite[Prop 2.44]{Joyce2013ACM} and Lemma \ref{lemembed}, it suffice to define gluing isomorphisms for $\mathbb{C}^\ast$-equivariant embeddings of charts of the form $(R,U,f,i)\hookrightarrow(R,U\times E,f\boxplus q,i\times 0)$, where $E$ is a finite dimensional vector space with linear $\mathbb{C}^\ast$ action and $q$ is a $\mathbb{C}^\ast$-invariant non-degenerate quadratic form. As seen above \eqref{diagnatident}, an orientation on $E$, which defines a descent isomorphism for $K^{1/2}_{X,s}$ for the embedding $(R,U,f,i)\hookrightarrow(R,U\times E,f\boxplus q,i\times 0)$, provides then naturally an orientation of $V^0$, which is a descent isomorphism for $K_{X^0,s^0}^{1/2}$ and the embedding $(R^0,U^0,f^0,i^0)\hookrightarrow(R^0,U^0\times E^0,f^0\boxplus q^0,i^0\times 0)$. The isomorphisms for $K_{X^0,s^0}^{1/2}$ glue on triple intersections if and only if those for $K_{X,s}^{1/2}$ glue also. $\Box$\medskip

Consider the decomposition of $X^0$ into connected components $X^0=\bigsqcup_{\pi\in\Pi}X^0_\pi$: one can further consider the oriented d-critical schemes $(X^0_\pi,s^0_\pi,K^{1/2}_{X^0_\pi,s^0_\pi})$. Denote by $\Ind_\pi$ the number of contracting weight in the tangent-obstruction complex of $X$ at $X^0_\pi$. In a critical chart $(R,U,f,i)$, the tangent-obstruction complex is given by $0\to TU\to T^\ast U\to0$, and then, denoting by $d_+,d_0,d_-$ respectively the number of contracting, invariant and repelling weights in $TU$ at $i(X^0_\pi\cap R)$, one has $\Ind_\pi=d_+-d_-$, \ie this definition is consistent with the previous definition on a critical chart.\medskip

\begin{theorem}\label{theospace}
    For $X$ an oriented d-critical algebraic space $X$ with an étale linearizable $\mathbb{C}^\ast$ action leaving the d-critical structure and the orientation invariant, one has natural isomorphisms of perverse sheaves and monodromic mixed Hodge modules on $X^0$ (with its orientation defined in Lemma \ref{lemorient}):
    \begin{align}
        \beta^\pm_{X,s}:(p_X^\pm)_!(\eta_X^\pm)^\ast P_{X,s}\overset{\simeq}{\to}\bigoplus_{\pi\in\Pi}\mathbb{L}^{\pm \Ind_\pi/2}P_{X^0_\pi,s^0_\pi}
    \end{align}
    which are compatible with polarization, in the sense that the following diagram is commutative:
    \[\begin{tikzcd}
        (p_X^+)_!(\eta_X^+)^\ast P_{X,s}\arrow[r,"\beta^+_{X,s}"]\arrow[d,"S_X^{-1}\circ (p_X^+)_!(\eta_X^+)^\ast(\Sigma_{X,s})"] & \bigoplus_{\pi\in\Pi}\mathbb{L}^{\Ind_\pi/2}P_{X^0_\pi,s^0_\pi}\arrow[d,"\bigoplus_{\pi\in\Pi}\mathbb{L}^{\Ind_\pi/2}\Sigma_{X^0_\pi,s^0_\pi}"]\\
        \mathbb{D}(p_X^-)_!(\eta_X^-)^\ast P_{X,s} & \mathbb{D}\bigoplus_{\pi\in\Pi}\mathbb{L}^{-\Ind_\pi/2}P_{X^0_\pi,s^0_\pi}\arrow[l,"\mathbb{D}\beta^-_{X,s}"]
    \end{tikzcd}\]
\end{theorem}

Proof: We can cover $X^0$ by critical charts of the form $(R^0,U^0,f^0,i^0)$, which are the $\mathbb{C}^\ast$-fixed part of $\mathbb{C}^\ast$-equivariant charts $(R,U,f,i)$, and the relate them by the transitive action of stabilization, \ie embedding of the form $\Phi^0:(R^0,U^0,f^0,i^0)\hookrightarrow(R^0,U^0\times E^0,f^0\boxplus q^0,i^0\times 0)$, which are the $\mathbb{C}^\ast$-fixed part of $\mathbb{C}^\ast$-equivariant embedding $\Phi:(R,U,f,i)\hookrightarrow(R,U\times E,f\boxplus q,i\times 0)$. Hence $\Perv(X^0)$ (or $MMHM(X^0)$) can be defined by descent datum with these charts and embeddings. Namely, $P_{X^0,s^0}$ is defined as said before by the descent datum $(\mathcal{PV}_{U,f},\Theta(\Phi^0))$. Because hyperbolic localization commutes with étale restriction, $(p_X^\pm)_!(\eta_X^\pm)^\ast P_{X,s}$ can be defined by the descent datum $((p_R^\pm)_!(\eta_R^\pm)^\ast\mathcal{PV}_{U,f},(p_R^\pm)_!(\eta_R^\pm)^\ast\Theta(\Phi))$.\medskip

On a chart $(R^0,U^0,f^0,i^0)$, we consider the isomorphism built in Proposition \ref{propchi}:
\begin{align}
        \beta^\pm_{U,f}:(p_R^\pm)_!(\eta_R^\pm)^\ast\mathcal{PV}_{U,f}\to\bigoplus_{\pi\in\Pi}\mathbb{L}^{\pm \Ind_\pi/2}\mathcal{PV}_{U^0_\pi,f^0_\pi}
\end{align}
We have to check he compatibility of these isomorphisms with stabilization. Consider as above $\Phi^0:(R^0,U^0,f^0,i^0)\hookrightarrow(R^0,U^0\times E^0,f^0\boxplus q^0,i^0\times 0)$, and consider the natural orientation on $E$ and $E^0$ defined by $K_{X,s}^{1/2}$ and $K_{X^0,s^0}^{1/2}$. consider now the following diagram:
\[\begin{tikzcd}[row sep=10ex,column sep=12ex]
    \scriptsize\parbox{6cm}{$(p_R^\pm)|_!(\eta_R^\pm)^\ast(\mathcal{PV}_{U\times E,f\boxplus q})$}\arrow[r,"\beta^\pm_{U\times E}"]\arrow[d,"\scriptsize\parbox{4cm}{$(p_R^\pm)|_!(\eta_R^\pm)^\ast\mathcal{TS}_{U,f,E,q}$}"] & 
    \scriptsize\parbox{6cm}{$\bigoplus_{\pi\in\Pi}\mathbb{L}^{\pm(\Ind_\pi}\mathcal{PV}_{U^0_{\pi}\times E^0,f^0_{\pi}\boxplus q^0}$}\arrow[d,swap,"\scriptsize\parbox{4cm}{$\bigoplus_{\pi\in\Pi}\mathbb{L}^{\pm(\Ind_\pi)/2}\mathcal{TS}_{U^0_{\pi},f^0_{\pi},E^0,q^0}$}"]\\
    \scriptsize\parbox{6cm}{$((p_R^\pm)_!(\eta_R^\pm)^\ast\mathcal{PV}_{U,f})\boxtimes(\mathcal{PV}_{E,q})$}\arrow[r,"\beta^\pm_U\times\beta^\pm_E"]\arrow[d,"\simeq"] &
    \scriptsize\parbox{6cm}{$\bigoplus_{\pi\in\mathbb{Z}}\mathbb{L}^{\pm \Ind_\pi/2}\mathcal{PV}_{U^0_{\pi},f^0_{\pi}})\boxtimes\mathcal{PV}_{E^0,q^0})$}\arrow[d,"\simeq"]\\
    (\eta_R^\pm)^\ast\mathcal{PV}_{U,f})\arrow[r,"\beta^\pm_U"] & \bigoplus_{\pi\in\mathbb{Z}}\mathbb{L}^{\pm \Ind_\pi/2}\mathcal{PV}_{U^0_{\pi},f^0_{\pi}})
    \end{tikzcd}\]
The above square is the commutative square of Proposition \ref{propcomthomseb} expressing the commutation of hyperbolic localization with Thom-Sebastiani isomorphism, and the square below is the tensor product of the commutative square \ref{diagnatident} with $(\eta_R^\pm)^\ast\mathcal{PV}_{U,f})$. Hence the above diagram is commutative, and is equal to:
\[\begin{tikzcd}[row sep=10ex,column sep=12ex]
    \scriptsize\parbox{6cm}{$(p_R^\pm)|_!(\eta_R^\pm)^\ast(\mathcal{PV}_{U\times E,f\boxplus q})$}\arrow[r,"\beta^\pm_{U\times E}"]\arrow[d,"\scriptsize\parbox{4cm}{$(p_R^\pm)|_!(\eta_R^\pm)^\ast\Theta(\Phi)$}"] & 
    \scriptsize\parbox{6cm}{$\bigoplus_{\pi\in\Pi}\mathbb{L}^{\pm(\Ind_\pi}\mathcal{PV}_{U^0_{\pi}\times E^0,f^0_{\pi}\boxplus q^0}$}\arrow[d,swap,"\scriptsize\parbox{4cm}{$\bigoplus_{\pi\in\Pi}\mathbb{L}^{\pm(\Ind_\pi)/2}\Theta(\Phi^0_\pi)$}"]\\
    (\eta_R^\pm)^\ast\mathcal{PV}_{U,f})\arrow[r,"\beta^\pm_U"] & \bigoplus_{\pi\in\mathbb{Z}}\mathbb{L}^{\pm \Ind_\pi/2}\mathcal{PV}_{U^0_{\pi},f^0_{\pi}})
\end{tikzcd}\]
Hence the $\beta^\pm_U$ are compatible with stabilization, and then by descent glue to give an isomorphism on $X^0$. for each $\mathbb{C}^\ast$-equivariant critical chart $(R,U,f,i)$ compatible with the monodromy and the duality. From Proposition \ref{propchi}, the $\beta^\pm_U$ are compatible with polarization and monodromy, hence by descent $\beta^\pm_X$ are also compatible (the commutativity of the corresponding squares can be checked locally). $\Box$\medskip

\begin{corollary}
    We define the cohomological DT invariant $[X]^{vir}:=[H_c(X,P_{X,s})]$ as the class of the hypercohomology with compact support of $P_{X,s}$ in the Grothendieck group of monodromic mixed Hodge structures. In the context of Theorem \ref{theospace}, one has:
\begin{align}
    [X]^{vir}=\sum_{\pi\in\Pi}\mathbb{L}^{\Ind_\pi/2}[X^0_\pi]^{vir}+[X-\eta(X^+)]^{vir}
\end{align}
\end{corollary}

Proof: $(X-\eta(X^+)\cup\bigcup_{\pi\in\Pi}X^+_\pi$ provides a stratification of $X$ by locally closed subvariety, hence one has:
\begin{align}
    [H_c(X,P_{X,s})]^{vir}&=[H_c(X^+,\eta^\ast P_{X,s})]+[H_c(X-\eta(X^+),P_{X,s})]\nn\\
    &=[H_c(X^0,p_!\eta^\ast P_{X,s})]+[H_c(X-\eta^\ast(X^+),P_{X,s})]\nn\\
    &=\bigoplus_{\pi\in\Pi}\mathbb{L}^{\Ind_\pi/2}[H_c(X^0,P_{X^0,s^0})]+[H_c(X-\eta^\ast(X^+),P_{X,s})]\nn\\
    &=\sum_{\pi\in\Pi}\mathbb{L}^{\Ind_\pi/2}[X^0_\pi]^{vir}+[X-\eta(X^+)]^{vir}
\end{align}
Here the first line holds because the class in the Grothendieck of monodromic mixed Hodge structures is a motivic invariant, the second line follows from the six functor formalism by viewing the cohomology with compact support as the proper pushforward to a point, the third line from the main result of this paper, and the last line by the definition of the cohomological DT invariants. $\Box$

\subsection{Compatibility with smooth pullbacks}\label{sectsmooth}

As mentioned in the subsection \ref{sectsemis}, the result of this subsection, the compatibility of the localization isomorphism with smooth pullbacks, is useful if one want to extends our result to DT theory on Artin stacks, following the approach of \cite{darbstack}. This was done in an earlier version of the paper, but removed here for clarity, but we have kept this subsection for future reference.\medskip

Consider a smooth map of d-critical algebraic spaces $\phi:(X,s)\to (Y,t)$ of relative dimension $d$ (\ie a smooth map such that $\phi^\star (t)=s$). \cite[Cor 3.8]{darbstack} shows that:
\begin{align}
    K_{X,s}=K_{Y,t}\otimes(\Lambda^{top}T^\ast_{X/Y})^{\otimes2}|_{X^{red}}
\end{align}
A smooth map of oriented d-critical scheme $\phi:(X,s,K_{X,s}^{1/2})\to (Y,t,K_{Yt}^{1/2})$ is then the data a smooth map $\phi:(X,s)\to (Y,t)$ of oriented d-critical scheme together with the data of an isomorphism:
\begin{align}\label{compatorient}
K^{1/2}_{X,s}=\phi^\ast(K^{1/2}_{Y,t})\otimes(\Lambda^{top}T^\ast_{X/Y})|_{X^{red}}. 
\end{align}
which is a square root of the above isomorphism. Given such a map, \cite[Prop 4.5]{darbstack} built an isomorphism:
\begin{align}
    \Delta_\phi:\phi^\ast[d]P_{Y,t}\simeq P_{X,s}
\end{align}

Consider now a $\mathbb{C}^\ast$-equivariant smooth map of oriented d-critical schemes $\phi:(X,s,K_{X,s}^{1/2})\to (Y,t,K_{Yt}^{1/2})$ of relative dimension $d$. From Lemma \ref{locmodsmooth}, $X$ is covered by  $\mathbb{C}^\ast$-equivariant critical charts $(R,U,f,i)$ such that there is a $\mathbb{C}^\ast$-equivariant critical chart and $(S,V,g,j)$ of $Y$ with a smooth map $\Phi:U\to V$ of relative dimension $d$ such that $f=g\circ\Phi$ and $\Phi\circ i=j\circ\phi$. Hence sheaves and perverse sheaves on $X^0$ can be defined by descent with a definition on each critical chart $(R^0,U^0,f^0,i^0)$ coming from a chart of this form, and comparison on charts coming from simultaneous stabilization of $(S,V,g,j)$ and $(R,U,f,i)$.

\begin{lemma}\label{lemorientsmooth}
    If $\phi:(X,s,K_{X,s}^{1/2})\to (Y,s,K_{Y,t}^{1/2})$ is a smooth $\mathbb{C}^\ast$ morphism of oriented d-critical loci, then $\phi^0:(X^0,s^0,K_{X^0,s^0}^{1/2})\to(Y^0,t^0,K_{Y^0,t^0}^{1/2})$ is aloso a smooth morphism of oriented d-critical loci.
\end{lemma}

Proof: The induced morphism $\phi^0:X^0\to Y^0$ is smooth of dimension $d$. One has:
\begin{align}
    (\phi^0)^\star(t^0)&=(\phi^0)^\star\xi_Y^\star(t)\nn\\
    &=\xi_X^\star\phi^\star(t)\nn\\
    &=\xi_X^\star(s)\nn\\
    &=s^0
\end{align}
where the first and the last line are the definition of $t^0$ and $s^0$, the second line follows from the functoriality of $^\star$ and the fact that $\xi_Y\phi^0=\phi\xi_X$, and the third line follows from the fact that $\phi$ is a morphism of d-critical structure. Then $\phi^0$ is a morphism of d-critical algebraic space.\medskip

We will now prove that:
\begin{align}\label{compatorientfix}
    K_{X^0,s^0}^{1/2}=(\phi^0)^\ast(K_{Y^0,t^0}^{1/2})\otimes(\Lambda^{top}T^\ast_{X^0/Y^0})|_{(X^0)^{red}}
\end{align}
We prove this equality by descent: this equality holds by the definitions $K_{X^0,s^0}^{1/2}|_{R^0}=K_{U^0}$ and $K_{Y^0,t^0}^{1/2}=K_{V^0}$ for each pair of critical charts of $X,Y$ of the form described above. Now, considering a stabilization for a pair of critical charts of $X,Y$, the compatibility of the descent isomorphisms on $X$ with \eqref{compatorient} implies the compatibility of the descent isomorphisms on $X^0$ with \eqref{compatorientfix}, which proves \eqref{compatorientfix}. $\Box$\medskip

 We can in particular consider the isomorphism:
\begin{align}
    \Delta_{\phi^0}:(\phi^0)^\ast[d]P_{Y^0,t^0}\to P_{X^0,s^0}
\end{align}

\begin{proposition}\label{compsmoothbeta}
    The hyperbolic localization isomorphism commutes with smooth pullbacks, namely for $\phi:(X,s,K_{X,s}^{1/2})\to (Y,t,K_{Y,t}^{1/2})$ a $\mathbb{C}^\ast$-equivariant smooth map of oriented critical algebraic space of relative dimension $d$, the following diagram of isomorphisms commutes:
    \[\begin{tikzcd}[row sep=5ex,column sep=20ex, scale cd=0.9]
    (\bigoplus_{\pi\in\Pi}\mathbb{L}^{\pm \Ind_\pi/2}(\phi^0_\pi)^\ast[d^0_\pi](p^\pm_Y)_!(\eta^\pm_Y)^\ast P_{Y,t})|_{R^0}\arrow[r,"\bigoplus_{\pi\in\Pi}\mathbb{L}^{\pm \Ind_\pi/2}(\phi^0_\pi)^\ast{[}d^0_\pi{]}\beta^\pm_{Y,t}"] & \bigoplus_{\pi,\pi'\in\Pi}\mathbb{L}^{\pm \Ind_{\pi'}/2}(\phi^0_\pi)^\ast[d^0_\pi]P_{Y^0_{\pi'-\pi},t^0_{\pi'-\pi}}\arrow[dd,swap,"\bigoplus_{\pi,\pi'\in\Pi}\mathbb{L}^{\pm \Ind_{\pi'}/2}\Delta_{\phi^0_{\pi}}"]\\
    ((p^\pm_X)_!(\eta^\pm_X)^\ast\phi^\ast[d]P_{Y,t}))|_{R^0}\arrow[d,"\Delta_\phi"]\arrow[u,"\simeq"] & \\
    ((p^\pm_X)_!(\eta^\pm_X)^\ast P_{X,s})|_{R^0}\arrow[r,"\beta^\pm_{X,s}"] & \bigoplus_{\pi'\in\Pi}\mathbb{L}^{\pm \Ind_{\pi'}/2}P_{X^0_{\pi'},s^0_{\pi'}}
    \end{tikzcd}\]
\end{proposition}

Proof: It suffices to show that this square commutes étale locally on $X^0$. According to lemma \ref{locmodsmooth}, we can then choose $\mathbb{C}^\ast$-equivariant critical charts $(R,U,f,i)$ and $(S,V,g,j)$ of $X$ and $Y$ containing $x$ and $y$ with a smooth map $\Phi:U\to V$ of relative dimension $d$ such that $f=g\circ\Phi$ and $\Phi\circ i=j\circ\phi$. On such charts, the commutation of the proposition is given by Proposition \ref{propcomsmooth}. $\Box$

\section{Technical lemmas}
\subsection{Trivialization of torus-equivariant embeddings}

\begin{lemma}\label{lemembed}
    Let $\Phi:(R,U,f,i)\to(S,V,g,j)$ be a $\mathbb{C}^\ast$-equivariant embedding of $\mathbb{C}^\ast$-equivariant critical charts. For $x\in i(R^0)$, there are smooth $\mathbb{C}^\ast$-equivariant algebraic space $U',V'$, $x'\in U'$, $\mathbb{C}^\ast$-equivariant morphisms $\iota: U'\to U$ with $\iota(x')=x$, $\jmath:V'\to V$, $\Phi':U'\to V'$, $\alpha:V'\to U$ and $\beta:V'\to E$ with $E$ a vector space with $\mathbb{C}^\ast$ action, with $\iota$, $\jmath$ and $\alpha\otimes\beta$ being étale, such that the following square commutes:

\[\begin{tikzcd}[row sep=10ex,column sep=20ex]
    U\arrow[d,"\Phi"] & U'\arrow[l,"\iota"]\arrow[r,"\iota"]\arrow[d,"\Phi'"] & U\arrow[d,"\Id_U\times 0"]\\V & V'\arrow[l,"\jmath"]\arrow[r,"\alpha\times\beta"] & U\times E
\end{tikzcd}\]
and $g\circ\jmath=f\circ q\circ\alpha:V'\to\mathbb{C}$, with $q$ a non-degenerate $\mathbb{C}^\ast$-invariant quadratic form on $E$.
\end{lemma}

Proof: We follow here closely the proof of \cite[Prop 2.22, 2.23]{Joyce2013ACM}, adapting it to the equivariant case: we take care to consider at each step only étale cover with $\mathbb{C}^\ast$ action, and to keep $\mathbb{C}^\ast$-equivariant functions and coordinates. The map $\Phi:U\to V$ is a $\mathbb{C}^\ast$-equivariant embedding of smooth schemes of respective dimension $\dim(U)=m,\dim(V)=m+n$. We can then find $\mathbb{C}^\ast$-equivariant étale coordinates $(\dot{y}_1,...,\dot{y}_m,\dot{z}_1,...,\dot{z}_n)$ on a $\mathbb{C}^\ast$-equivariant Zariski open neighborhood $\dot{V}$ of $j(x)$ in $V$ such that $j(x) = (0,..., 0)$ and $\Phi(U)\cap\dot{V}$ is the locus $\dot{z}_1=...=\dot{z}_n=0$ in $\dot{V}$. Set $\dot{U}=\Phi^{-1}(\dot{V})$ and $\dot{x}_a=\dot{y}_a\circ\Phi|_{\dot{U}}$ for $a = 1, . . . , m$. Then $\dot{U}$ is an open neighborhood of $i(x)$ in $U$ , and $(\dot{x}_1, ...,\dot{x}_m)$ are étale coordinates on $\dot{U}$ with $i(x) = (0,...,0)$. Then the ideal $I_{R,U}=I_{(df)}$ is on $\dot{U}$ the ideal generated by the $\mathbb{C}^\ast$-equivariant functions $\frac{\partial f}{\partial\dot{x}_a}$ for $a = 1, . . . , m$, and the ideal $I_{S,V}= I_{(dg)}$ is on $\dot{V}$ the ideal generated by the $\mathbb{C}^\ast$-equivariant functions $\frac{\partial g}{\partial\dot{y}_a}$ for $a = 1, . . . , m$ and $\frac{\partial g}{\partial\dot{z}_b}$ for $b = 1, . . . , n$. Since $\Phi$ maps $U$ to $\dot{z}_1=...=\dot{z}_n=0$ and $i(R)$ to $j(R)\subset j(S)$, we have $I_{(df)}\simeq I_{(dg)} |_{\dot{z}_1=...=\dot{z}_n=0}$ , that is:
\begin{align}
    (\frac{\partial f}{\partial\dot{x}_a}(\dot{y}_1,...,\dot{y}_m):a=1,...,m)=(&\frac{\partial g}{\partial\dot{y}_a}(\dot{y}_1,...,\dot{y}_m,0,...,0):a=1,...,m,\nn\\&\frac{\partial g}{\partial\dot{z}_b}(\dot{y}_1,...,\dot{y}_m,0,...,0):b=1,...,n)
\end{align}
this holds provided each $\frac{\partial g}{\partial\dot{z}_b}(\dot{y}_1,...,\dot{y}_m,0,...,0)$ lies in $(\frac{\partial g}{\partial\dot{y}_a}(\dot{y}_1,...,\dot{y}_m,0,...,0):a=1,...,m$. Thus, making $\dot{U}$, $\dot{V}$ smaller if necessary, we can suppose there exist étale locally $\mathbb{C}^\ast$-equivariant functions $A_{ab}(\dot{y}_1,...,\dot{y}_m)$ on $\dot{U}$ for $a = 1, . . . , m, b = 1, . . . , n$ such that for each b:
\begin{align}
    \frac{\partial g}{\partial\dot{z}_b}(\dot{y}_1,...,\dot{y}_m,0,...,0)=\sum_{a=1}^m A_{ab}(\dot{y}_1,...,\dot{y}_m)\frac{\partial g}{\partial\dot{y}_a}(\dot{y}_1,...,\dot{y}_m,0,...,0)
\end{align}

Defining the $\mathbb{C}^\ast$-equivariant functions $\tilde{y}_a=\dot{y}_a-\sum_{a=1}^m A_{ab}(\dot{y}_1,...,\dot{y}_m)\dot{z}_b$, and $\tilde{z}_b=\dot{z}_b$, they give $\mathbb{C}^\ast$-equivariant étale coordinates on a $\mathbb{C}^\ast$-invariant neighborhood $\tilde{V}$ of $j(x)$. We also define the $\mathbb{C}^\ast$-invariant subspace $\tilde{U}=\Phi^{-1}(\tilde{V})$, and define on it the $\mathbb{C}^\ast$-equivariant étale coordinates $\tilde{x}=\dot{x}|_{\tilde{U}}$: then $\tilde{y}_a\circ\Phi|_{\tilde{U}}=\tilde{x}_a$ and $\tilde{z}_b\circ\Phi|_{\tilde{U}}=0$. Then:
\begin{align}
    \frac{\partial g}{\partial\tilde{z}_b}(\tilde{y}_1,...,\tilde{y}_m,0,...,0)=&\sum_{a=1}^m\frac{\partial g}{\partial\dot{y}_a}(\dot{y}_1,...,\dot{y}_m,0,...,0).\frac{\partial\dot{y}_a}{\partial\tilde{z}_b}+\sum_{c=1}^n\frac{\partial g}{\partial\dot{z}_c}(\dot{y}_1,...,\dot{y}_m,0,...,0).\frac{\partial\dot{z}_c}{\partial\tilde{z}_b}\nn\\
   =&\sum_{a=1}^m\frac{\partial g}{\partial\dot{y}_a}(\dot{y}_1,...,\dot{y}_m,0,...,0)(-A_{ab}(\dot{y}_1,...,\dot{y}_m))\nn\\&+\sum_{c=1}^n(\sum_{a=1}^nA_{ac}(\dot{y}_1,...,\dot{y}_m)).\frac{\partial g}{\partial\dot{y}_a}(\dot{y}_1,...,\dot{y}_m,0,...,0).\delta_{bc}\nn\\
    =& 0
\end{align}
We define the $\mathbb{C}^\ast$-equivariant space $\check{V}$ and the $\mathbb{C}^\ast$-equivariant étale morphisms $\check{\jmath},\check{\alpha},\check{\beta}$ by the $\mathbb{C}^\ast$-equivariant Cartesian square:
\[\begin{tikzcd}[row sep=10ex,column sep=20ex]
    \check{V}\arrow[r,"\check{\jmath}"]\arrow[d,"\check{\alpha}\times\check{\beta}"] & \tilde{V}\arrow[d,"\scriptsize\parbox{2cm}{$(\tilde{y}_1,...,\tilde{y}_m,\tilde{z}_1,...,\tilde{z}_n)$}"]\\
    \tilde{U}\times\mathbb{C}^n \arrow[r,"\scriptsize\parbox{2cm}{$(\tilde{x}_1,...,\tilde{x}_m)\times\Id_{\mathbb{C}^n}$}"] & \mathbb{C}^{m+n}
\end{tikzcd}\]
where $\mathbb{C}^{m+n}$ is provided with the $\mathbb{C}^\ast$-action induced from the action of $\mathbb{C}^\ast$ on the étale coordinates. There is a unique $\check{v}\in\check{V}$ with $\check{\alpha}(\check{v})=i(x)$, $\check{\beta}(\check{v})=(0,...,0)$ and $\check{\jmath}(\check{v})=j(x)$.). We can regard $\check{\jmath}:\check{V}\to\tilde{V}\subset V$ as an étale $\mathbb{C}^\ast$-invariant open set in V. Define $\mathbb{C}^\ast$-equivariant étale coordinates $(\check{y}_1,...,\check{y}_m,\check{z}_1,...,\check{z}_n):\check{V}\to\mathbb{C}^{m+n}$ by $\check{y}_a=\tilde{y}_a\circ\check{\jmath},\check{z}_b=\tilde{z}_b\circ\check{\jmath}$, and define the $\mathbb{C}^\ast$-invariant functions $\check{f},\check{g}:\check{V}\to\mathbb{C}$ by $\check{f}=f\circ\check{\alpha},\check{g}=g\circ\check{\jmath}$, and $\check{h}=\check{g}-\check{f}:\check{V}\to\mathbb{C}$. The previous argument now shows that on the smooth subscheme $\check{U}\subset\check{V}$ defined by $\check{z}_1=...=\check{z}_n=0$ we have $\check{h}|_{\check{U}}$ and $\frac{\partial\check h}{\partial\check{z}_b}|_{\check{U}}=0$ for $b=1,...,n$. Therefore the $\mathbb{C}^\ast$ invariant function $\check{h}$ lies in the ideal $(\check{z}_1,...,\check{z}_n)^2$ generated by $\mathbb{C}^\ast$-equivariant functions on $\check{V}$. So making $\check{U},\check{V},\tilde{U},\tilde{V}$ smaller but still $\mathbb{C}^\ast$-invariant, we may write $\check{h}=\sum_{b,c=1}^n\check{z}_b\check{z}_c Q_{bc}$ for some $\mathbb{C}^\ast$-invariant $Q_{bc}:\check{V}\to\mathbb{C}$ with $Q_{bc}=Q_{cb}$. Up to a $\mathbb{C}^\ast$-equivariant linear change of coordinates $\check{z}$, we can write:\medskip
\begin{align}
    \sum_{b,c=1}^n\check{z}_b\check{z}_c Q_{bc}(0,...,0)=2\sum_{b=1}^r\check{z}_{2b-1}\check{z}_{2b}+\sum_{b=2r+1}^n\check{z}_b^2
\end{align}
with $n_{2i-1}=-n_{2i}>0$ for $i\leq r$ and $n_i=0$ for $i\geq 2r+1$. Suppose that $n\geq 2r+1$: the $\mathbb{C}^\ast$-invariant function $Q_{nn}(\check{y}_1,...\check{y}_m,\check{z}_1,...,\check{z}_n)$ is invertible near $\check{v}$, we can take a double étale cover of $\check{V}-Q_{nn}^{-1}(0)$ such that $Q_{nn}$ has a square root $Q_{nn}^{1/2}$. Because $Q_{nn}$ is $\mathbb{C}^\ast$-invariant, $\check{V}-Q_{nn}^{-1}(0)$ is stable under $\mathbb{C}^\ast$, hence we can extend the $\mathbb{C}^\ast$ action to the double cover, such that the cover is $\mathbb{C}^\ast$-equivariant. We replace then $\check{V}$ by this étale neighboorhood of $\check{v}$, and all the coordinates and $Q_{bc}$, and $Q_{nn}^{1/2}$, are still $\mathbb{C}^\ast$-equivariant. We can then write:
\begin{align}
    h&=\sum_{b,c=1}^n\check{z}_b\check{z}_c Q_{bc}\nn\\
   &= \sum_{b,c=1}^{n-1}\check{z}_b\check{z}_c(Q_{bc}-Q_{nn}^{-1}Q_{bn}Q_{cn})+(Q_{nn}^{1/2}\check{z}_n+\sum_{b=1}^{n-1}Q_{nn}^{-1/2}Q_{bn}\check{z}_b)^2\nn\\
   &=\sum_{b,c=1}^{n-1}\check{z}_b\check{z}_c\hat{Q}_{bc}+z_n^2
\end{align}
where $\hat{Q}_{bc}=Q_{bc}-Q_{nn}^{-1}Q_{bn}Q_{cn}$ is $\mathbb{C}^\ast$ equivariant of character $-n_b-n_c$ and $z_n=Q_{nn}^{1/2}\check{z}_n+\sum_{b=1}^{n-1}Q_{nn}^{-1/2}Q_{bn}\check{z}_b$ is $\mathbb{C}^\ast$-invariant, with $\frac{\partial z_n}{\partial\check{z}_b}(\check{v})=\delta_{nb}$ and $\hat{Q}_{bc}(0,...,0)=\delta_{bc}$. By recursion we can find an étale open neighborhood of $\check{v}$ in $V$, $\mathbb{C}^\ast$-invariant functions on it $z_b$ for $b\geq 2r+1$ with $\frac{\partial z_b}{\partial\check{z}_c}(\check{v})=\delta_{bc}$, and $\mathbb{C}^\ast$-equivariant functions on it $\check{Q}_{bc}$ for $b,c\geq 2r$ with $\check{Q}_{bc}(0,...,0)=Q_{bc}(0,...,0)$, such that:
\begin{align}
    h=\sum_{b,c=1}^{2r}\check{z}_b\check{z}_c\check{Q}_{bc}+\sum_{b=2r+1}^n z_b^2
\end{align}
On a $\mathbb{C}^\ast$-invariant open neighborhood of $\check{v}$ such that $\check{Q}_{2r-1,2r}$ is invertible:
\begin{align}
    h=&\sum_{b,c=1}^{2r-2}\check{z}_b\check{z}_c(\check{Q}_{bc}-\sum_{b,c=1}^{2r-2}\check{Q}_{2r-1,2r}^{-1}\check{Q}_{b,2r}\check{Q}_{2r-1,c})+\sum_{b=2r+1}^n z_b^2\nn\\&+2(\check{z}_{2r-1}+\sum_{b=1}^{2r-2}\check{Q}_{2r-1,2r}^{-1}\check{Q}_{b,2r}\check{z}_b)(\check{Q}_{2r-1,2r}\check{z}_{2r}+\sum_{c=1}^{2r-2}\check{Q}_{2r-1,c}\check{z}_c)\nn\\
    =&\sum_{b,c=1}^{2r-2}\check{z}_b\check{z}_c\tilde{Q}_{bc}+2z_{2r-1}z_{2r}+\sum_{b=2r+1}^n z_b^2
\end{align}
where $\tilde{Q}_{bc}=\check{Q}_{bc}-\sum_{b,c=1}^{2r-2}\check{Q}_{2r-1,2r}^{-1}\check{Q}_{b,2r}\check{Q}_{2r-1,c}$ is $\mathbb{C}^\ast$-equivariant of character $-n_b-n_c$ with $\tilde{Q}_{bc}(0,...,0)=Q_{bc}(0,...,0)$, $z_{2r-1}=\check{z}_{2r-1}+\sum_{b=1}^{2r-2}\check{Q}_{2r-1,2r}^{-1}\check{Q}_{b,2r}\check{z}_b$ is $\mathbb{C}^\ast$ equivariant of character $n_{2r-1}$ with $\frac{\partial z_{2r-1}}{\partial\check{z}_b}(\check{v})=\delta_{2r-1,b}$, and $z_{2r}=\check{Q}_{2r-1,2r}\check{z}_{2r}+\sum_{c=1}^{2r-2}\check{Q}_{2r-1,c}\check{z}_c$ is $\mathbb{C}^\ast$-equivariant of weight $n_{2r}$ with $\frac{\partial z_{2r}}{\partial\check{z}_b}(\check{v})=\delta_{2r,b}$. By recursion, there is then on an a $\mathbb{C}^\ast$-invariant étale open neighborhood $\jmath':V'\to\check{V}$ of $\check{v}$ and $\mathbb{C}^\ast$-equivariant functions $z_b$ with $\frac{\partial z_b}{\partial\check{z}_c}(\check{v})=\delta_{bc}$ (\ie by restricting $V'$ to a $\mathbb{C}^\ast$-invariant open neighborhood of $\check{v}$, $(y_1:=y_1\circ\jmath',...,y_m:=y_m\circ\jmath',z_1,...,z_n$ forms a system of étale coordinates of $V'$), such that:
\begin{align}\label{eqh}
    h=2\sum_{b=1}^rz_{2b-1}z_{2b}+\sum_{b=2r+1}^nz_b^2
\end{align}
Define the $\mathbb{C}^\ast$-invariant subset $U'=\{v'\in V':z_1(v')=...=z_n(v')=0$. Define the $\mathbb{C}^\ast$-equivariant maps $\iota:U'\to U,\jmath:V'\to V,\Phi':U'\to V',\alpha:V'\to U,\beta:V'\to\mathbb{C}^n$ by $\iota:=\check{\alpha}\circ\jmath'|_U'$, $\jmath:=\check{\jmath}\circ\jmath'$, $\Phi'=\Id_{U'}$, $\alpha=\check{\alpha}\circ\jmath'$ and $\beta=(z_1,...,n_n)$. As $\jmath':V'\to\check{V}$
is an étale open neighborhood of $\check{v}$ in $\check{V}$ , there exists $u'\in V'$ with $\jmath'(u')=\check{v}$,
and $z_b(u'=0$ for $b=1,...,n$ as $\check{z}_b(\check{v})=0$, so $u'\in U$
with $\iota(u')= \check{\alpha}\circ\jmath'(u') = \check{\alpha}(\check{v}) = i(x)$. Also $\iota,\jmath,\alpha\times\beta$ are étale as $\check{alpha}\times\check{\beta},\check{\jmath},\jmath'$ are. The same computations as in the proof of \cite[Prop 2.23]{Joyce2013ACM} show that in a neighborhood of $u'$, which can be taken to be $\mathbb{C}^\ast$-invariant, we have $\Phi\circ\iota=\jmath\circ\Phi$, so this making $u'$ smaller and still $\mathbb{C}^\ast$-invariant it holds on $U'$. The
equations $\alpha\circ\Phi'=\iota,\beta\circ\Phi'=0$ are immediate, and $g\circ\jmath= f\circ\alpha+(2\sum_{b=1}^rz_{2b-1}z_{2b}+\sum_{b=2r+1}^nz_b^2)\circ\beta$ follows from \eqref{eqh}. $\Box$

\subsection{Smooth torus-equivariant morphism of d-critical scheme}

\begin{lemma}\label{locmodsmooth}
    Consider a smooth morphism of d-critical scheme $\phi:(X,s)\to (Y,t)$ of relative dimension $d$. For $x\in X$,  we can then choose two $\mathbb{C}^\ast$-equivariant critical charts $(R,U,f,i)$ and $(S,V,g,j)$ of $X$ and $Y$ containing $x$ and $y:=f(x)$ with a smooth map $\Phi:U\to V$ of relative dimension $d$ such that $f=g\circ\Phi$ and $\Phi\circ i=j\circ\phi$.
\end{lemma}

Proof: This is a $\mathbb{C}^\ast$-equivariant version of the result used in the proof of \cite[Prop 4.5]{darbstack}: as before, we check that we can do each step in a $\mathbb{C}^\ast$-equivariant way. As in the proof of \cite[Prop 2.43]{Joyce2013ACM}, because the $\mathbb{C}^\ast$ action on $Y$ is étale-locally linearizable, we can find an étale neighborhood $S'$ of $y$ with a closed $\mathbb{C}^\ast$-equivariant embedding into a smooth algebraic space with $\mathbb{C}^\ast$-action $V'$, with $\dim(V')=\dim(T_y Y)$. According to the proof of the proposition \ref{prophyplocsmoothaff}, the $\mathbb{C}^\ast$-equivariant smooth map $\phi^{-1}(S')\to S'$ is étale locally a $\mathbb{C}^\ast$-equivariant affine fibration, hence we can find a $\mathbb{C}^\ast$-equivariant étale neighborhood $S$ of $y$, $\mathbb{C}^\ast$-equivariant embeddings $i:R:=\phi^{-1}(S)\to U$, $j:S\to V$ into smooth algebraic spaces with $\mathbb{C}^\ast$-action $U,V$ with dimension $\dim(U)=\dim(T_x X)$, $\dim(V)=\dim(T_y Y)$, with a $\mathbb{C}^\ast$-equivariant smooth map of relative dimension $d$ $\Phi:U\to V$ such that $\Phi\circ i=j\circ\phi|_R$. Hence as in the proof of \cite[Prop 2.43]{Joyce2013ACM}, up to shrinking $S,V$ (and then $R,U$) to a $\mathbb{C}^\ast$-invariant open subspace we can find a $\mathbb{C}^\ast$ invariant regular function $g:V\to\mathbb{C}$ such that $\iota_{S,V}(t|_S)=j^{-1}(g)+I_{S,V}^2$, and $S=\Crit(g)$, \ie $(S,V,g,j)$ is a $\mathbb{C}^\ast$-equivariant critical chart of $Y$ near $y$. Defining the $\mathbb{C}^\ast$-equivariant regular function $f:=g\circ\Phi:U\to\mathbb{C}$, because $s=\phi^\star(t)$, one has $\iota_{R,U}(s|_R)=i^{-1}(f)+I_{R,U}^2$. Also $\Phi\circ i=j\circ\phi|_R$ and $\Phi,\phi$ are smooth of relative
dimension $d$, then $i(R)=\Phi^{-1}(j(S))$ in a neighborhood of $i(x)$: taking the union of all these neighborhood, one can shrink $U$ (and then $R$) to $\mathbb{C}^\ast$-invariant open neighborhood satisfying $i(R)=\Phi^{-1}(j(S))$. Because $\Phi$ is smooth, one has:
\begin{align}
    \Crit(f)=\Phi^{-1}(\Crit(g))=\Phi^{-1}(j(S))=i(R)
\end{align}
\ie $(R,U,f,i)$ is a $\mathbb{C}^\ast$-equivariant critical chart near $x$. We have then build $\mathbb{C}^\ast$-equivariant critical charts $(R,U,f,i)$ and $(S,V,g,j)$ of $X$ and $Y$ near $x$ and $y$ with a smooth map $\Phi:U\to V$ of relative dimension $d$ such that $f=g\circ\Phi$ and $\Phi\circ i=j\circ\phi$. $\Box$

\bibliography{ref}
\bibliographystyle{hamsplain}
\end{document}